\documentclass[reqno]{amsart}

\numberwithin{equation}{section}
\newtheorem{Theorem}{Theorem}[section]
\newtheorem{Lemma}[Theorem]{Lemma}
\newtheorem{Proposition}[Theorem]{Proposition}

\newcommand{\thref}[1]{Theorem \ref{#1}}
\newcommand{\leref}[1]{Lemma \ref{#1}}
\newcommand{\prref}[1]{Proposition \ref{#1}}
\newcommand{\reref}[1]{Remark \ref{#1}}

\newcommand{\boldZ}{{\mathbb{Z}}}

\theoremstyle{definition}
\newtheorem{Remark}[Theorem]{Remark}

\newcommand{\Span}{{\mathrm{span}}}                              %
\newcommand{\C}{\mathbb{C}}                                      %
\newcommand{\R}{\mathbb{R}}                                      %
\newcommand{\G}{{\Gamma}}                                        %
\newcommand{\Gt}{{\tilde{\Gamma}}}                               %
\newcommand{\Phit}{\tilde{\Phi}}                                   %
\newcommand{\cK}{\mathcal{K}}                                    %
\newcommand{\cL}{\mathcal{L}}                                    %
\newcommand{\cM}{\mathcal{M}}                                    %
\newcommand{\nn}{\nonumber}
\def\TT{{\mathbb T}}                                      %
\begin{document}
\title[Two Variable orthogonal polynomials on the bi-circle]{Two variable
orthogonal polynomials on the bi-circle and structured matrices}

\author[J.~S.~Geronimo]{Jeffrey~S.~Geronimo}\thanks{The authors were partially
supported by an NSF grant and a NATO collaborative linkage grant}
\address{School of Mathematics, Georgia Institute of Technology,
Atlanta, GA 30332--0160, USA}
\email{geronimo@math.gatech.edu}

\author[H.~Woerdeman]{Hugo~Woerdeman}
\address{Department of Mathematics, Drexel University,
Philadelphia, PA 19104, USA} \email{Hugo.Woerdeman@drexel.edu}
\maketitle
\date{March 26, 2006}

\begin{abstract} We consider bivariate polynomials orthogonal on the bicircle
with respect to a positive linear functional. The lexicographical
and reverse lexicographical orderings are used to order the
monomials. Recurrence formulas are derived between the polynomials
of different degrees. These formulas link the orthogonal
polynomials constructed using the lexicographical ordering with
those constructed using the reverse lexicographical ordering.
Relations between the coefficients in the recurrence formulas are
derived and used to give necessary and sufficient conditions for
the existence of a positive linear functional. These results are
then used to  construct a class of two variable measures supported
on the bicircle that are given by one over the magnitude squared
of a stable polynomial. Applications to Fej\'er-Riesz
factorization are also given.
\end{abstract}

\bigskip

\noindent {\bf Keywords}: Bivariate orthogonal polynomials,
positive definite linear functionals, moment problem, doubly
Toeplitz matrices, recurrence coefficients.

\bigskip

\noindent {\bf MSC}: 42C05 , 30E05 , 47A57, 15A48, 47B35

\section{Introduction}                                        %

Bivariate polynomials orthogonal on the bicircle have been
investigated mostly in the electrical engineering community in
relation to the design of stable recursive filters for
two-dimensional filtering. In particular we note the work of Genin
and Kamp \cite{GK} who were interested in the following problem.
Given any two variable polynomial $q(z,w)$ with $q(0,0)\ne0$, let
$a_{k,l}(z,w)$ be its planar least squares inverse polynomial of
degree $(k,l)$; i.e., $a_{k,l}$ minimizes the mean quadratic value
of $1-a_{k,l}q$ on the bicircle. What properties does $a_{k,l}$
have? At the time it was conjectured  the minimizing polynomials
were stable, i.e. $a_{k.l}(z,w)\ne0,\ |z|<1,\ |w|<1$, which they
showed was false. Their investigation was carried further by
Delsarte, Genin and Kamp \cite{DGKls} who developed the connection
between these polynomials and matrix polynomials orthogonal on the
unit circle \cite{DGKm}. In the development of this connection
these authors were lead to examine moment matrices that were block
Toeplitz matrices where each block entry is itself a Toeplitz
matrix. Such structured matrices are called doubly Toeplitz
matrices and arise naturally in the bivariate trigonometric moment
problem. These types of matrices arose more recently in the work
of Geronimo and Woerdeman \cite{GW1} in their investigation of the
bivariate Fej\'er-Riesz factorization theorem. These authors were
able to resolve the question when a strictly positive bivariate
trigonometric polynomial of a certain degree can be written as the
magnitude squared of a stable polynomial of the same degree. In
this work the authors used the fact that the theory of orthogonal
polynomials on the unit circle provides a proof of the one
variable Fej\'er-Riesz theorem which does not use the fundamental
theorem of algebra. We intend here to continue to investigate the
properties of bivariate polynomials orthogonal on the bicircle and
clarify their role in the Fej\'er-Riesz Theorem.

A major difficulty encountered in the theory of orthogonal
polynomials of more than one variable is which monomial ordering
to use. For bivariate real orthogonal  polynomials the preferred
ordering is the total degree ordering which is the one set by
Jackson \cite{J}. For polynomials with the same total degree the ordering
is lexicographical. As noted in  Delgado et al \cite{DGIM} in
their study of orthogonal polynomials associated with doubly
Hankel matrices there is a good reason for choosing this ordering
which is that if new orthogonal polynomials of higher degree are
to be constructed then their orthogonality relations will not
affect the relations governing the lower degree polynomials.
However in order for the moment matrix to be doubly Toeplitz the
monomial orderings that need to be used are lexicographical and
reverse lexicographical.

We begin in Section 2 by considering finite dimensional
subspaces spanned by the monomials $z^i w^j, |i|\le n, \, |j|\le
m$, and exhibiting the connection between positive linear
functionals defined on this space and positive definite doubly
Toeplitz matrices. We then introduce certain matrix orthogonal
polynomials and show how they give the Cholesky factors for the
inverse of the doubly Toeplitz matrices considered above. The
results in \cite{GW1} show that these polynomials play a role in
the parametric moment problem. In Section 3 we construct two
variable orthogonal polynomials, where the monomials are ordered
according to the lexicographical ordering. When these polynomials
are organized into vector orthogonal polynomials they can be
related to the matrix orthogonal polynomials constructed
previously. From this relation it is shown that these vector
polynomials are the minimizers of a certain quadratic functional.
Using the orthogonality relation, recurrence relations satisfied
by the vector polynomials and their counterparts in the reverse
lexicographical ordering are derived and relations between these
recurrence coefficients are exhibited. In Section 4 a number of
Christoffel-Darboux like formulas are derived.  In Section 5 we
use the relations between the coefficients derived in Section 3 to
develop an algorithm to construct the coefficients in the
recurrence formulas at a particular level $(n,m)$ say, in terms of
the coefficients at the previous levels plus a certain number of
unknowns. The collection of these unknowns is in one to one
correspondence with the number of moments needed to construct the
vector polynomials up to level $(n,m)$. This is used in Section 6
to construct a positive linear functional from the recurrence
coefficients. The construction allows us to find necessary and
sufficient conditions on the recurrence coefficients for the
existence of a positive linear functional which is in one to one
correspondence with the set of positive definite doubly Toeplitz
matrices. In Section 7 we examine conditions under which the
linear functional can be represented as a positive measure
supported on the bicircle having the form of one over the
magnitude squared of a stable polynomial. This gives a new proof of
the Fej\'er-Riesz result of \cite{GW1}. Finally in Section 8
examples are given that illustrate various aspects of the theory
developed.

\bigskip
\section{Positive linear functionals and Doubly Toeplitz matrices}     %
In this section we consider moment matrices associated with the
lexicographical ordering which is defined by
$$
(k,\ell)<_{\rm lex} (k_1,\ell_1)\Leftrightarrow k<k_1\mbox{ or }
(k=k_1\mbox{ and } \ell<\ell_1),
$$
and the reverse lexicographical ordering defined by
$$
(k,\ell)<_{\rm revlex} (k_1,\ell_1)\Leftrightarrow
(\ell,k)<_{\rm lex} (\ell_1,k_1).
$$
Both of these orderings are linear orders and in addition they satisfy
$$
(k,\ell)<(m,n)\Rightarrow (k+p,\ell+q)<(m+p,n+q).
$$
In such a case, one may associate a halfspace with the ordering
 which is defined by
 $\{ (k,l) \ : \  (0,0)  <  (k,l) \}$. In the case of the
 lexicographical ordering we shall denote the associated halfspace by $H$
 and refer to it
 as {\it the standard halfspace}. In the case of the reverse lexicographical
 ordering we shall denote the associated halfspace by $\tilde H$.
 Instead of starting with the ordering, one may also start with a halfspace
 $\hat H$ of $\boldZ^2$ (i.e., a set $\hat H$
 satisfying $\hat H+\hat H \subset \hat H$, $\hat H\cap (-\hat H)= \emptyset$,
 $\hat H \cup (-\hat H) \cup \{ (0,0)\} = \boldZ^2$) and define an ordering
 via
 $$ (k,l) <_{\hat H} (k_1,l_1) \Longleftrightarrow (k_1-k,l_1-l) \in \hat H . $$
 We shall refer to the order $<_{\hat H}$ as
 {\it the order associated with} $\hat H$. Note that the lexicographical and reverse
lexicographical orderings do not respect total degree.

Let $\prod^{n,m}$ denote the bivariate Laurent linear subspace
$\Span \{z^iw^j,\, -n\le i\le n,\,-m\le j\le m\}$.  Let
$\cL_{n,m}$ be a linear functional defined on $\prod^{n,m}$ by
$$
\cL_{n,m}(z^{-i}w^{-j})=c_{i,j}= \overline{\cL(z^i w^j)}.
$$
We will call $c_{i,j}$ the $(i,j)$ moment of $\cL_{n,m}$ and
$\cL_{n,m}$ a moment functional.  If we form the
$(n+1)(m+1)\times(n+1)(m+1)$ matrix $C_{n,m}$ for $\cL_{n,m}$ in
the lexicographical ordering then, as noted in the introduction,
it has the special block Toeplitz form
\begin{equation}\label{toepone}
C_{n,m} = \left[
\begin{matrix}
C_0 & C_{-1} & \cdots & C_{-n}
\\
C_1 & C_0 & \cdots & C_{-n+1}
\\
\vdots &  & \ddots & \vdots \cr C_n & C_{n-1} & \cdots & C_{0}
\end{matrix}
\right],
\end{equation}
where each $C_i$ is an $(m+1)\times(m+1)$ Toeplitz matrix as
follows
\begin{equation}\label{toeptwo}
C_i=\left[
\begin{matrix}
c_{i,0}& c_{i,-1} & \cdots & c_{i,-m}
\\
\vdots & &\ddots & \vdots \cr c_{i,m} &  & \cdots &
c_{i,0}
\end{matrix}
\right],\qquad i=-n,\dots, n.
\end{equation}
Thus $C_{n,m}$ has a doubly Toeplitz structure. If the reverse
lexicographical ordering is used in place of the lexicographical
ordering we obtain another moment matrix $\tilde C_{n,m}$ where
the roles of $n$ and $m$ are interchanged.

Let us introduce the notion of centro-transpose symmetry. We
denote the transpose of a matrix $A$ by $A^T$. A square matrix $A$
is said to be {\it centro-transpose symmetric} if $JAJ=A^T$ where
$J$ is the matrix with 1's on the antidiagonal and zeros
elsewhere. Note that a Toeplitz matrix is centro-transpose
symmetric. We have the following useful lemmas which characterize
Toeplitz and doubly Toeplitz matrices in terms of centro-transpose
symmetry.

\begin{Lemma}\label{toep1}
An $(n+1)\times (n+1)$ matrix $A=(a_{i,j})_{i,j=0}^{n}$ is
Toeplitz if and only if both $A$ and $\hat A:=(a_{i,j})_{i,j=0}^{n-1}$ are centro-transpose symmetric.
\end{Lemma}

\begin{proof}
Notice that $JAJ=A^T$ is equivalent to $a_{n-i,n-j} = a_{j,i}$,
$0\le i,j \le n$. Similarly, the centro-transpose symmetry of
$\hat A$ is equivalent to $a_{n-1-i,n-1-j}=a_{j,i}$, $0\le i,j\le
n-1$. But then
$$ a_{i+1,j+1} = a_{n-j-1,n-i-1}= a_{i,j} , 0\le i,j\le n-1 , $$
and thus it follows that $A$ is Toeplitz.

As $A$ and $\hat A$ are Toeplitz, the converse is immediate.
\end{proof}

\begin{Lemma}\label{toepblock}
Let $A=(A_{i,j})$, $i,j=1,\dots, k$, where each $A_{i,j}$ is a
complex $m\times m$ matrix. Then $A$ is a doubly Toeplitz matrix
if and only if $A^T=JAJ$, $A_1^T=J_1 A_1 J_1$ and $A_2^T=J_1 A_2
J_1$. Here $A_1$ is obtained from $A$ by deleting the last block
row and column and $A_2$ is obtained from $A$ by removing the last
row and column of each $A_{i,j}$. The  matrices $J$
and $J_1$ are square matrices of appropriate size with ones on the
antidiagonal and zeros everywhere else.
\end{Lemma}

\begin{proof}
Again the necessary conditions follow from the structure of $A$.
To see the converse note that $A^T=JAJ$ implies that $A_{j,i}^T =J_2 A_{k-i,k-j}J_2$, where $J_2$ is the $m\times m$ matrix with ones on the reverse diagonal and zeros everywhere else . This coupled with the condition on $A_1$ implies that $A$ is a block Toeplitz matrix from \leref{toep1} and $J_2 A_{i,j} J_2= A_{i,j}^T$. These relations plus the condition on $A_2$ and \leref{toep1} gives the result.
\end{proof}

\begin{Remark} The conclusion of the above lemmas hold if  we replace deleting the last (last block)  row and column by deleting the first (first block) row and column.\end{Remark}

We say that the moment functional $\cL_{n,m}:\prod^{n,m}\to\C$
is positive definite or positive semidefinite if
\begin{equation}\label{pfunct}
\cL_{n,m} (|p|^2)>0\quad {\rm or} \quad\cL_{n,m}(|p|^2)\ge 0
\end{equation}
for every every nonzero polynomial $p\in\prod^{n,m}$. It
follows from a simple quadratic form argument that $\cL_{n,m}$ is
positive definite or positive semidefinite if and only if its
moment matrix $C_{n,m}$ is positive definite or positive
semidefinite respectively.

We will say that $\cL$ is positive definite or positive
semidefinite if
$$
\cL(|p|^2)>0\qquad {\rm or}\qquad\cL (|p|^2)\ge 0
$$
respectively for all nonzero polynomials.  Again these conditions
are equivalent to the moment matrices $C_{n,m}$ being positive
definite or positive semidefinite for all positive integers $n$
and $m$. The above discussion leads to,

\begin{Lemma}\label{btoeptwo}
Let $C_{n,m}$ be a positive (positive semi-) definite
$(n+1)(m+1)\times(n+1)(m+1)$ matrix given by \eqref{toepone} and
\eqref{toeptwo}. Then there is a positive (positive semi-) definite moment
functional $\cL_{n,m}:\prod^{n,m}\to\C$ associated with
$C_{n,m}$ given by
$$
c_{i,j} = \cL_{n,m}(z^{-i}w^{-j})=\overline{\cL_{n,m}(z^i w^j)},\qquad -n\le i\le n,\qquad -m\le j\le
m.
$$
The converse also holds.
\end{Lemma}

Let $\prod^{n}_{m+1}$ be the set of all $(m+1)\times(m+1)$ complex
valued matrix polynomials of degree $n$ or less, $\prod_{m+1}$ the set of all $(m+1)\times(m+1)$ complex
valued matrix polynomials, and $M^{m,n}$ is
the space of $m\times n$ matrices. For a matrix $M$ we let
$M^\dagger$ denote the conjugate transpose (or the adjoint) of
$M$. For a polynomial $Q(z,w)$ we let $Q^\dagger (z,w)$ denote the
polynomial in $z^{-1}$ and $w^{-1}$ defined by $Q^\dagger (z,w) =
Q(\frac{1}{\overline z}, \frac{1}{\overline w} )^\dagger$. If the
positive moment functional $\cL_{n,m}:\prod^{n,m}\to\C$ is
extended to two variable polynomials with matrix coefficients in
the obvious way, we can associate to it a positive matrix function
$\cL_m: \prod^{n}_{m+1}\times\prod^{n}_{m+1}\to M^{m+1,m+1}$
defined by,
\begin{equation}\label{relfunct}
\left[\cL_m(P(z),Q(z))\right]_{i,j}=\cL_{n,m}(\left[P(z,w)\ Q^{\dagger}(z,w)\right]_{i,j}),\quad 1\le i,j\le m+1
\end{equation}
where
$$
P(z,w)=P(z)
\left[\begin{matrix} 1\\ \vdots\\
w^m\end{matrix}\right]\  {\rm and}\
Q(z,w)=Q(z)\left[\begin{matrix} 1\\ \vdots\\
w^m\end{matrix}\right].
$$

Equation~\eqref{relfunct} shows that if $\cL_{n,m}$ can be represented as a positive measure $\mu$ supported on the bicircle then for $f$ an $m+1\times m+1$ matrix function continuous on the unit circle,
$$
\cL_m(f)=\int_{-\pi}^{\pi}f(\theta)dM_m(\theta),
$$
where $M_m$ is the $m+1\times m+1$ matrix measure given by,
$$
dM_m(\theta)=\int_{\phi=-\pi}^{\pi}\left[\begin{matrix} 1\\ \vdots\\
w^m\end{matrix}\right]d\mu(\theta,\phi)\left[\begin{matrix} 1\\ \vdots\\
w^m\end{matrix}\right]^{\dagger},
$$
which shows that $M_m$ is Toeplitz.

Because of the structure of $C_{n,m}$ we can associate to $\cL_m$
matrix valued orthogonal polynomials in the following manner
\cite{DGKm}, \cite{DGKls}, \cite{GW1}. Let $\{R^m_i(z)\}^n_{i=0}$
and $\{L^m_i(z)\}^n_{i=0}$ be $(m+1)\times(m+1)$ complex valued
matrix polynomials given by
\begin{equation}\label{rightpoly}
R^m_i(z)=R^m_{i,i}z^i+R^m_{i,i-1}z^{i-1}+\cdots, \qquad i=0,\dots,
n,
\end{equation}
and
\begin{equation}\label{leftpoly}
L^m_i(z)=L^m_{i,i}z^i+L^m_{i,i-1}z^{i-1}+\cdots, \qquad i=0,\dots,
n,
\end{equation}
satisfying
\begin{equation}\label{rightorth}
\cL_m ({R^m_i}^{\dagger}, {R^m_j}^{\dagger})=\delta_{ij}I_{m+1}
\end{equation}
and
\begin{equation}\label{leftorth}
\cL_m ({L^m_i}, {L^m_j})=\delta_{ij}I_{m+1}
\end{equation}
respectively, where $I_{m+1}$ denotes the $(m+1)\times(m+1)$
identity matrix. The above relations uniquely determine the
sequences $\{R^m_i\}^n_{i=0}$ and $\{L^m_i\}^n_{i=0}$ up to a unitary
factor and this factor will be fixed by requiring $R^m_{i,i}$ and $L^m_{i,i}$ to be upper
triangular matrices with positive diagonal entries. We write
\begin{equation}\label{Lpoly}
L^m_i (z)=[0\cdots\ 0\ L^m_{i,i}\ L^m_{i,i-1}\ \cdots\ L^m_{i,0}]
\left[\begin{matrix} z^n I_{m+1}\\ z^{n-1}I_{m+1}\\ \vdots\\ I_{m+1}
\end{matrix}\right]
\end{equation}
and
\begin{equation}\label{Lpolyone}
\hat L^m_n(z)=\left[\begin{matrix}L^m_n(z)\\ L^m_{n-1}(z)\\ \vdots\\ L^m_0(z)
\end{matrix}\right] =L\left[\begin{matrix}z^n I_{m+1}\\ z^{n-1}I_{m+1}\\
\vdots\\ I_{m+1}\end{matrix}\right],
\end{equation}
where
\begin{equation}\label{Lpolytwo}
L=\left[\begin{matrix}
L^m_{n,n} & L^m_{n,n-1} & \cdots  & L^m_{n,0}\\
0 & L^m_{n-1,n-1} & \cdots & L^m_{n-1,0}\\
\vdots & & \ddots \\
0 & 0& \cdots & L^m_{0,0}\end{matrix}\right].
\end{equation}
In an analogous fashion write,
\begin{equation}\label{Rpolyone}
\hat R^m_n(z)=\left[\begin{matrix}R^m_0(z)\\ R^m_1(z)\\ \vdots\\ R^m_n(z)
\end{matrix}\right] =\left[\begin{matrix}I_{m+1}&\ldots& z^nI_{m+1}\end{matrix}\right]R,
\end{equation}
where
\begin{equation}\label{Rpolytwo}
R=\left[\begin{matrix}
R^m_{0,0} & R^m_{1,0}& \cdots  & R^m_{n,0}\\
0 & R^m_{1,1} & \cdots & R^m_{n,1}\\
\vdots & & \ddots \\
0 &0 & \cdots & R^m_{n,n}\end{matrix}\right].
\end{equation}
By lower (respectively upper) Cholesky factor $A$ (respectively $B$) of a
positive definite matrix $M$ we mean
\begin{equation}\label{chol}
M=AA^{\dagger} = BB^{\dagger},
\end{equation}
where $A$ is a lower triangular matrix with positive diagonal
elements and $B$ is an upper triangular matrix with positive
diagonal elements. With the above we have the following well known lemma \cite{KVM},

\begin{Lemma}\label{lemchol}
Let $C_{n,m}$ be a positive definite block Toeplitz matrix given by \eqref{toepone}
then $L^{\dagger}$ is the lower Cholesky factor and $R$ is the upper Cholesky factor
of $C^{-1}_{n,m}$.
\end{Lemma}

\begin{proof}
To obtain \eqref{chol} note that \eqref{leftorth} implies that
$$
I=\cL_m((\hat L^m_n), (\hat L^m_n))=L\cL_m\left(\left[
\begin{matrix}
z^n I_{m+1}\\ z^{n-1}I_{m+1}\\ \vdots \\ I_{m+1}
\end{matrix}\right],
\left[\begin{matrix} z^n I_{m+1}\\ z^{n-1}I_{m+1}\\ \vdots \\
I_{m+1}\end{matrix}\right]\right)L^{\dagger}=LC_{n,m}L^{\dagger},
$$
where $I$ is the $(n+1)(m+1)\times(n+1)(m+1)$ identity matrix.
Since $C_{n,m}$ is invertible we find,
$$
C^{-1}_{n,m} = L^{\dagger}L.
$$
The result for $R$ follows in an analogous manner.
\end{proof}
  From this formula and \eqref{Lpolytwo} we find,
\begin{equation}\label{L_n}
L^m_n(z) = \left[({L^m_{n,n}}^{\dagger})^{-1},0,0,\dots 0\right]
C^{-1}_{n,m}
[z^n I_{m+1},z^{n-1}I_{m+1},\ldots, I_{m+1}]^T,
\end{equation}
and
\begin{equation}\label{R_n}
R^m_n(z) = \left[I_{m+1},zI_{m+1},\ldots, z^nI_{m+1}\right]
C^{-1}_{n,m}
\left [0, 0, \dots, 0,  ({\bar R^m_{n,n}})^{-1}\right]^T.
\end{equation}
Note that ${L^m_{n,n}}^{\dagger}$ is the lower Cholesky factor of
$[I_{m+1},0,\cdots,0]C^{-1}_{n,m}\ [I_{m+1},0,\cdots,0]^T$ while
$R^m_{n,n}$ is the upper Cholesky factor of
$[0,\cdots,I_{m+1}]C^{-1}_{n,m}\ [0,\cdots,I_{m+1}]^T$.

The theory of matrix orthogonal polynomials (\cite{DGKm},
\cite{KVM}, \cite{Ro}, \cite{Simon}) can be applied to obtain the recurrence
formulas
\begin{equation}\label{mrecurrence}
\begin{split}
A_{i+1,m} L^m_{i+1}(z) &= z L^m_{i}(z) - E_{i+1,m}{\overleftarrow R}^m_{i}(z)\quad i= 0,\ldots, n-1,
\\
R^m_{i+1}(z)\hat A_{i+1,m} &= z R^m_{i}(z) - {\overleftarrow L}^m_{i}(z)E_{i+1,m}\quad i= 0,\ldots, n-1,
\end{split}
\end{equation}
where
\begin{equation}\label{eis}
E_{i+1,m}=\cL_{m}(z L^m_i, \overleftarrow {R}^m_{i})=\cL_{m}{(\overleftarrow{L}^m_i}^{\dagger} , \left(z R^m_i)^{\dagger}\right)
\end{equation}
and
\begin{equation}\label{bis}
\begin{split}
A_{i+1,m}=\cL_{m}(zL^m_i,\ L^m_{i+1})=L^m_{i,i}(L^m_{i+1,i+1})^{-1}\\
\hat A_{i+1,m}= \cL_{m}(R^{m\dagger}_{i+1},
\left(z R^m_i)^{\dagger}\right)=(R^m_{i+1,i+1})^{-1} R^m_{i,i}.
\end{split}
\end{equation}
For a matrix polynomial $B$ of degree $n$ in $z$, $\overleftarrow B(z)=z^n B^\dagger(1/z)$.
By multiplying the first equation in \eqref{mrecurrence} on the left by $\bar z {L^m_i(z)}^{\dagger}$ and the second equation on the right by $\bar z {R^m_i(z)}^{\dagger}$ then integrating we see that
\begin{equation}\label{firstk}
\begin{split}
A_{i+1,m}A_{i+1,m}^{\dagger}= I_{m+1}-E_{i+1,m}E_{i+1,m}^{\dagger},\\
{\hat A_{i+1,m}}^{\dagger}\hat A_{i+1,m} = I_{m+1}-E_{i+1,m}^{\dagger} E_{i+1,m}.
\end{split}
\end{equation}
The above equations and the properties of $A_{i+1,m}$ and $\hat
A_{i+1,m}$ show that $E_{i+1,m}$ is a strictly contractive matrix
and that $A_{i+1,m}$ is the upper Cholesky factor of
$I_{m+1}-E_{i+1,m}E_{i+1,m}^{\dagger}$. Similarly $\hat
A_{i,m}^{\dagger}$ is the lower Cholesky factor of
$I_{m+1}-E_{i+1,m}^{\dagger} E_{i+1,m}$. Furthermore \eqref{bis}
and \eqref{firstk} show that
\begin{equation}\label{determ}
\det((L^m_{i+1,i+1})^{\dagger}L^m_{i+1,i+1})^{-1}=\det(C_0)\prod_{j=1}^{i+1}\det(I_{m+1}-E_{j,m} E^{\dagger}_{j,m}).
\end{equation}
The recurrence formulas \eqref{mrecurrence} can be inverted in the following manner. Multiply the reverse of the second equation in \eqref{mrecurrence} on the right by $E_{i+1,m}$ to obtain
$$
E_{i+1,m}\hat A_{i+1,m}^{\dagger}\overleftarrow{R}^m_{i+1}(z) =
E_{i+1,m}\overleftarrow{R}^m_{i}(z)
-z E_{i+1,m}E^{\dagger}_{i+1,m} L^m_{i}(z).
$$
Add this equation to the first equation in \eqref{mrecurrence} then use \eqref{firstk} to eliminate $A_{i+1,m}$ and $\hat A_{i+1,m}^{\dagger}$ to find
\begin{equation}\label{frinv}
(A_{i+1,m}^{\dagger})^{-1}L^m_{i+1}(z)+ E_{i+1,m}(\hat A_{i+1,m})^{-1}\overleftarrow{R}^m_{i+1}(z)=z L^m_{i}(z).
\end{equation}
In a similar manner we find
\begin{equation}\label{secinv}
R^m_{i+1}(z)(\hat A_{i+1,m}^{\dagger})^{-1}+ \overleftarrow{L}^m_{i+1}(z)(A_{i+1,m})^{-1} E_{i+1,m}=z R^m_{i}(z).
\end{equation}
From the recurrence formulas it is not difficult to derive the
Christoffel-Darboux formulas \cite{DGKm},
\begin{equation}\label{mCDeq}
\begin{split}
{\overleftarrow{R}^m_k(z)}^{\dagger}{\overleftarrow{R}^m_k(z_1)}-\bar z z_1{L^m_k(z)}^{\dagger} L^m_k(z_1)=(1-\bar z z_1)\sum_{i=0}^k{L^m_i(z)}^{\dagger} L^m_i(z_1)\\
{\overleftarrow{L}^m_k(z_1)}{\overleftarrow{L}^m_k(z)}^{\dagger}-\bar z z_1 R^m_k(z_1){R^m_k(z)}^{\dagger}=(1- \bar z z_1)\sum_{i=0}^k R^m_i(z_1){R^m_i(z)}^{\dagger}
\end{split}
\end{equation}

These formulas give rise to the matrix Gohberg-Semencul formulas
\cite{GH},\cite{KVM} when the linear equations obtained by
equating like powers of ${\bar z}^i z^j$ are put in matrix form.
Some properties that follow from the above formulas \cite{DGKm}
[Theorems 9, 14 and 15] are that $\overleftarrow{R}^m_i(z)$ and
$\overleftarrow{L}^m_k(z)$ have empty kernels for $|z|\le 1$,
i.e.,
\begin{equation}\label{stablem}
\det(\overleftarrow{R}^m_i(z))\ne0\ne\det(\overleftarrow{L}^m_k(z)),\
|z|\le1.
\end{equation}
Such polynomials are called stable matrix polynomials and if we
write
\begin{equation}\label{wone}
W_k(z)=
\left[{\overleftarrow{L}^m_k(z)}{\overleftarrow{L}^m_k(z)}^{\dagger}\right]^{-1},
\end{equation}
and
$$
C^k_j=\frac{1}{2\pi}\int_{-\pi}^{\pi} e^{-ij\theta}
W_k(e^{i\theta})d\theta
$$
then
\begin{equation}\label{fcoef}
C^k_j= C_j, \quad  |j|\le k.
\end{equation}
Furthermore
\begin{equation}\label{wt2}
W_k=\left[{\overleftarrow{R}^m_k(z)}^{\dagger}{\overleftarrow{R}^m_k(z)}\right]^{-1}.
\end{equation}
If $\overleftarrow{L}^m_k(z)$ $(\overleftarrow{R}^m_k(z))$
satisfies \eqref{stablem} and \eqref{fcoef} we will say it is
stable and has spectral matching (up to level $k$).  Another
useful result shown in \cite{DGKm} is
\begin{equation}\label{entropy}
\log\det((L^m_{i+1,i+1})^{\dagger}L^m_{i+1,i+1})^{-1}=\frac{1}{2\pi}\int_{-\pi}^{\pi}
\log \det W_k(\theta) d\theta .
\end{equation}
From the stability of $\overleftarrow{R}^n_{i+1}$ and
$\overleftarrow{L}^m_{i+1}$, \eqref{frinv} and \eqref{secinv} give
the following formulas for the recurrence coefficients
$E_{i+1,m}$:
\begin{align}
E_{i+1,m}&=-(A_{i+1,m}^{\dagger})^{-1}L^m_{i+1}(0)\overleftarrow{R}^m_{i+1}(0)^{-1}\hat
A_{i+1,m}\\ \nn
&=-A_{i+1,m}\overleftarrow{L}^m_{i+1}(0)^{-1}R^m_{i+1}(0)(\hat
A_{i+1,m}^{\dagger})^{-1} . \label{reflect}\end{align} We also
note that $\overleftarrow{L}^m_k(z)$ and
$\overleftarrow{R}^m_k(z)$ are minimizers of certain quadratic
functions. To see this denote the set of $(m+1)\times (m+1)$
hermitian matrices as ${\rm Herm}(m+1)$ and let
$\cM : \prod_{m+1} \to {\rm Herm}(m+1)$ be given by,
\begin{equation}\label{lfunt}
\cM[X(z)]={\cL}_m(X, X)-(X(0)+X(0)^{\dagger}),
\end{equation}
then Delsarte, Genin and Kamp have shown \cite{DGKm} that for a
given degree $k$, ${\cM}$ is minimized by
$\overleftarrow{L}^m_k(z)L^m_{k,m}$ with value
$(L^m_{k,m})^{\dagger}L^m_{k,m}$. Likewise  $\hat{\cM} :
\prod_{m+1} \to {\rm Herm}(m+1)$
given by,
\begin{equation}\label{rfunt}
\hat{\cM}[X(z)]={\cL}_m(X^{\dagger}, X^{\dagger})-(X(0)+X(0)^{\dagger}),
\end{equation}
is minimized by  $R^m_{k,m}\overleftarrow{R}^m_k(z)$ and takes the value $R^m_{k,m}
(R^m_{k,m})^{\dagger}$. Thus we find,
\begin{equation}\label{decrl}
(L^m_{k,m})^{\dagger}L^m_{k,m}\ge (L^m_{k+1,m})^{\dagger}L^m_{k+1,m}
\end{equation}
and
\begin{equation}\label{decrr}
R^m_{k,m}(R^m_{k,m})^{\dagger}\ge
R^m_{k+1,m}(R^m_{k+1,m})^{\dagger} .
\end{equation}
Here $A\ge B$ for two $(m+1)\times(m+1)$ matrices means that $A-B$
is positive semidefinite. The above discussion leads to Burg's
entropy Theorem. Consider the class of $M^m$ of  $(m+1)\times
(m+1)$ matrix Borel  measures on the unit circle and for each such
measure $\mu$ write the Lebesgue decomposition of $\mu=
\mu_{ac}+\mu_s$ where $d\mu_{ac}/d\theta = W(\theta)$. Let $S^m_n$
be the subset of $M^m$ such that each $\mu\in S^m_n$ has the same
Fourier coefficients $C_i,\ |i|\le n$ and ${\mathcal
E}(\mu)=\frac{1}{2\pi}\int_{-\pi}^{\pi} \ln\det(W) d\theta
>-\infty$. Then there is a unique measure which maximizes the
above entropy function ${\mathcal E}(\mu)$ and this measure is
given by $d\mu=W(\theta)d\theta$ with
$W(\theta)=Q^m_n(\theta)^{-1}$ where $Q^m_n(\theta)$ is a positive
$(m+1)\times (m+1)$ matrix trigonometric polynomial of degree $n$.

This leads to a simple proof of the Matrix Fej\'er-Reisz
factorization Theorem (Helson \cite{H}, Dritschel \cite{D},
McLean-Woerdeman \cite{MW}, Geronimo-Lai \cite{GL}) which will be
useful later.
\begin{Lemma}\label{mfr}
Let $Q^m_n(\theta)$ be a strictly positive $(m+1)\times (m+1)$
matrix trigonometric polynomial then
$Q^m_n(\theta)=\overleftarrow{L}_n^m(z)
(\overleftarrow{L}^m_n(z))^{\dagger},\ z=e^{i\theta}$ where
$\overleftarrow{L}_n^m$ is a stable $(m+1)\times (m+1)$ matrix
polynomial of degree $n$. Furthermore $L^m_n$ is given by
\eqref{L_n}.
\end{Lemma}
\begin{proof}  Since $Q^m_n(\theta)$ is strictly positive we can compute
the moments \hfill\break
$C_j=\frac{1}{2\pi}\int_{-\pi}^{\pi} e^{-ij\theta} Q^m_n(\theta)^{-1}d\theta$. If we compute
the matrix orthogonal polynomials associated with
these Fourier coefficients we find that $W_n$ has spectral matching up to $n$. That is its
Fourier coefficients match $C_i$ for $|i|\le n$. The Maximum Entropy Theorem implies that
$ Q^m_n(\theta)= W_n^{-1}$, which gives the result.
\end{proof}

The matrix Fej\'er-Riesz Theorem now follows.
\begin{Theorem}\label{matrixfr}
Let $Q^m_n(\theta)\ge0$ be a positive $(m+1)\times (m+1)$ matrix
trigonometric polynomial then $Q^m_n(\theta)=P_n^m(z)
(P^m_n(z))^{\dagger},\ z=e^{i\theta}$ where $P_n^m$ is an outer
(nonzero for $|z|<1$) $(m+1)\times (m+1)$ matrix polynomial.
\end{Theorem}
\begin{proof} Let $Q^m_{n,\epsilon} = \epsilon I + Q^m_n\, \epsilon>0$
then $Q^m_{n,\epsilon}$ satisfies the hypotheses of the above lemma.
Thus  $Q^m_{n,\epsilon}=P^m_{n,\epsilon}(P^m_{n,\epsilon})^{\dagger}$. The proof now follows
by taking the limit as $\epsilon$ tends to zero.
\end{proof}

It was observed by Delsarte et. al. \cite{DGKls} that if the $C_k$ in $C_{n,m}$ are
centro-transpose symmetric then
\begin{equation}\label{hicoeff}
({L^m_{i,i}}^{\dagger}L^m_i(z))^T=J_m R^m_{i}(z){R^m_{i_i}}^{\dagger}J_m\quad
 i=0,\ldots,n
\end{equation}
where $J_m$ is the $(m+1)\times (m+1)$  matrix with ones on the
reverse diagonal and zeros everywhere else. This can easily be
seen from formulas \eqref{L_n} and \eqref{R_n} since in this case
from \leref{toepblock} $C^T_{n,m}=J C_{n,m}J$ with $J$ the
$(n+1)(m+1)\times (m+1)(n+1)$ matrix with ones down the
anti-diagonal and zeros everywhere else . This leads to the
following characterization of positive definite doubly Toeplitz
matrices in terms of certain recurrence coefficients. We will
denote by $C^m_0$ the $m\times m$ matrix obtained from $C_0$ by
eliminating the first row and first column of $C_0$.

\begin{Theorem}\label{reflcentro} Suppose $C_{n,m}$ is positive definite.
Then the Fourier coefficients $C_i, |i|\le n$, are
centro-transpose symmetric if and only if $E_{k,m}, k=1, \ldots,
n$, and $C_0$ are centro-transpose symmetric. Consequently,
$C_{n,m}$ is doubly Toeplitz if and only if $E_{k,i}, k=1, \ldots
, n, i = m-1,m$, $C_0$ and $C_0^m$    are centro-transpose
symmetric.
\end{Theorem}

\begin{proof} Examining the leading coefficients in \eqref{hicoeff} and using the
fact that $L^m_{i,i}$ and $R^m_{i,i}$ are upper triangular we find
that (see also \cite{DGKls}), $(L^m_{i,i})^T=J_m R^m_{i,i} J_m$
for $i=0,\dots , n$. Thus
\begin{equation}\label{jsym}
L^m_i(z)^T=J_m R^m_i(z)J_m,\ i=0,\ldots,n
\end{equation}
From \eqref{eis} we find, since $J_m J_m =I_{m+1}$,
$$
J_m E_{i+1,m}J_m=\cL_{m}(\bar z J_m{L^m_i}^{\dagger}J_m,J_m{\overleftarrow {R}^m_{i}}^{\dagger}J_m)=\cL_{m}(\overleftarrow{L}^m_i,z R^m_{i})^T=E_{i+1,m}^T
$$
To show the converse note that if $E_{i,m}$ is centro-transpose
symmetric then from \eqref{firstk} we obtain $$
J_m(A_{i,m}A_{i,m}^{\dagger})^T J_m=J_m(I_m-
E_{i,m}E_{i,m}^{\dagger})J_m =\overline{(I- E_{i,m}^{\dagger}
E_{i,m})}={\hat A_{i,m}}^T\overline{\hat A_{i,m}},
$$
so that 
\begin{equation}\label{centroahat}
J_m A_{i,m}J_m=\hat A_{i,m}^T.
\end{equation}
Since $C_0$ is
centro-transpose symmetric and $L_{0,m}(z)$ is the upper Cholesky
factor of $C_0$, we see
 that $J_m L^m_0 J_m={R^m_0}^T$.
Thus by induction using \eqref{mrecurrence} we find that
$J_m L^m_n(z)J_m=R^m_n(z)^\top$. The first part of the result now follows from the spectral
matching of $W_n$ (\eqref{fcoef}) and \eqref{wt2} . The second part of the Theorem follows by
applying the above argument to $C_{n,m-1}$ and
$C^m_0$ then
using \leref{toep1} .
\end{proof}

In the next two sections we present recurrence formulas and an algorithm that computes recurrence
coefficients for a positive definite doubly Toeplitz matrix.

\section{Bivariate orthogonal polynomials}
In this section we examine the properties of two variable
orthogonal polynomials where the monomial ordering is either
lexicographical or reverse lexicographical. The study of
orthogonal  polynomials on the bicircle with this ordering was
begun by Delsarte et.al. \cite{DGKls} and extended in \cite{GW1}.
Given a positive definite linear functional $\cL_{N,M}
:\prod^{N,M}\to\C$ we perform the Gram-Schmidt procedure using the
lexicographical ordering and define the orthonormal polynomials
$\phi_{n,m}^l(z,w),\ 0\le n\le N,\, 0\le m\le M, \, 0\le l\le m,$
by the equations,
\begin{equation}\label{sorthogonal}
\begin{split}
&\cL_{N,M}(\phi_{n,m}^l z^{-i}w^{-j})=0, \quad  0\le i<n\ {\rm and }\
0\le j\le m\ {\rm\ or}\ i=n\ {\rm and }\ 0\le j< l,
\\
&\cL_{N,M}(\phi_{n,m}^l(\phi_{n,m}^l)^{\dagger})=1,
\end{split}
\end{equation}
and
\begin{equation}\label{sorthogonaldeg}
\phi_{n,m}^{l}(z,w) = k^{n,l}_{n,m,l} z^n w^l + \sum_{(i,j)<_{\rm
lex}(n,l)} k^{i,j}_{n,m,l}z^iw^j.
\end{equation}
With the convention $k^{n,l}_{n,m,l}>0$, the above equations uniquely
specify $\phi^l_{n,m}$. Polynomials orthonormal with respect to
$\cL_{N,M}$ but using the reverse lexicographical ordering will be
denoted by $\tilde \phi^l_{n,m}$. They are uniquely determined by the
above relations with the roles of $n$ and $m$ interchanged.

Set
\begin{equation}\label{vectpoly}
\Phi_{n,m}=\left[\begin{matrix} \phi_{n,m}^{m}\\ \phi_{n,m}^{m-1}\\
\vdots\\ \phi_{n,m}^{0} \end{matrix}\right] =
K_{n,m}\left[\begin{matrix} z^n w^m\\ z^n w^{m-1}\\ \vdots\\ 1\end{matrix}\right],
\end{equation}
where the $(m+1)\times(n+1)(m+1)$ matrix $K_{n,m}$ is given by
\begin{equation}\label{knmmatrix}
K_{n,m}=\left[\begin{matrix} k_{n,m,m}^{n,m}&
k_{n,m,m}^{n,m-1}&\cdots & \cdots&\cdots& k_{n,m,m}^{0,0} \\ 0&
k_{n,m,m-1}^{n,m-1}&\cdots & \cdots&\cdots&
k_{n,m,m-1}^{n,0}\\\vdots &\ddots&\ddots&\ddots&\ddots&\ddots
\\0&\cdots& k_{n,m,0}^{n,0}& k_{n,m,0}^{n-1,m}&\cdots &
k_{n,m,0}^{0,0}\end{matrix}\right] .
\end{equation}
As indicated above denote
\begin{equation}\label{vectpolytilde}
\tilde\Phi_{n,m} =\left[\begin{matrix} \tilde \phi_{n,m}^{n}\\ \tilde
\phi_{n,m}^{n-1}\\ \vdots\\ \tilde \phi_{n,m}^{0}\end{matrix}\right]
=\tilde K_{n,m}\left[\begin{matrix}w^m z^n \\  w^m z^{n-1}\\ \vdots\\ 1
\end{matrix}\right],
\end{equation}
where the $(n+1)\times(n+1)(m+1)$ matrix $\tilde K_{n,m}$ is
given similarly to \eqref{knmmatrix} with the roles of $n$ and $m$
interchanged. For the bivariate polynomials $\phi^l_{n,m}(z,w)$ above
we define the reverse polynomials $\overleftarrow \phi^l_{n,m}(z,w)$ by the relation
\begin{equation}\label{reversepoly}
\overleftarrow \phi^l_{n,m}(z,w)=z^n w^m\bar{\phi}_{n,m}^l(1/z,1/w).
\end{equation}
With this definition $\overleftarrow \phi^l_{n,m}(z,w)$ is again a
polynomial in $z$ and $w$, and furthermore
\begin{equation}\label{reversevect}
\overleftarrow\Phi_{n,m}(z,w):= z^n w^m \Phi_{n,m}^{\dagger}(1/z,1/w)=\left[\begin{matrix} \overleftarrow\phi_{n,m}^{m}\\ \overleftarrow\phi_{n,m}^{m-1}\\
\vdots\\ \overleftarrow\phi_{n,m}^{0} \end{matrix}\right].
\end{equation}
An analogous procedure is used to define $\overleftarrow{\tilde\phi}^l_{n,m}$.

In order to ease the notation to find recurrence formulas for the vector
polynomials $\Phi_{n,m}$ we introduce the inner product,
\begin{equation}\label{innerprod}
\langle X, Y\rangle=\cL_{N,M}(X Y^{\dagger}).
\end{equation}

Let $\hat\prod^{n,m}$ be the linear span of $z^i w^j,\ 0\le i\le n,0\le j\le m$, $\hat\prod_k^{n,m}$ be
the vector space of $k$ dimensional
vectors with entries in $\hat\prod^{n,m}$, and $\hat\prod^{m}_{m+1}=\hat\prod^{\infty,m}_{m+1}$.

Utilizing the orthogonality relations \eqref{sorthogonal} we
obtain the following auxiliary results.

\begin{Lemma}\label{lemvectorth}
Suppose $\Phi\in\hat\prod_{k}^{n,m}$. If
$\Phi$ satisfies the orthogonality relations,
\begin{equation}\label{vectorth}
\langle\Phi, z^i w^j\rangle = 0,\quad 0\le i<n,\quad 0\le j\le m,
\end{equation}
then $\Phi = T\Phi_{n,m}$, where $T$ is a $k\times(m+1)$ matrix.
If $k=m+1$, $T$ is upper triangular with positive diagonal
entries, and if $\langle\Phi,\Phi\rangle=I_{m+1}$, then
$T=I_{m+1}$.
\end{Lemma}

\begin{Lemma}\label{lemvectorthtilde}
Suppose $\Phit\in\hat\prod_k^{n,m}$. If
$\Phit$ satisfies the orthogonality relations,
\begin{equation}\label{vectorthtilde}
\langle\Phit, z^i w^j\rangle = 0,\quad 0\le i\le n,\quad 0\le j< m,
\end{equation}
then $\Phit = T\Phit_{n,m}$, where $T$ is an $k\times(n+1)$
matrix. If $k=n+1$, $T$ is upper triangular with positive diagonal
entries, and if $\langle\Phit,\Phit\rangle=I_{n+1}$, then
$T=I_{n+1}$.

\end{Lemma}

With the above we can make contact with the matrix orthogonal
polynomials introduced in Section 2. This was observed by Delsarte
et. al. \cite{DGKls}

\begin{Lemma}\label{lemrelmatrix}
Let $\Phi_{n,m}$ be given by \eqref{vectpoly}. Then
\begin{equation}\label{relmatrix}
\Phi_{n,m}= L^m_n(z)[w^m,w^{m-1},\ldots,1]^T,
\end{equation}
\begin{equation}\label{revmatrix}
\overleftarrow{\Phi}_{n,m}= [1,w,\ldots,w^m]J_m\overleftarrow{R}^m_n(z)^TJ_m,
\end{equation}
and
\begin{align}
\nonumber\left[\begin{matrix}\Phi_{n,m}(z,w)\\ \Phi_{n-1,m}(z,w)\\ \vdots\\
\Phi_{0,m}(z,w)\end{matrix}\right]& =\left[\begin{matrix}L^m_n(z)\\
L^m_{n-1}(z)\\ \vdots\\ L^m_0(z)\end{matrix}\right] [w^m,w^{m-1},\ldots,
1]^T\\& =L\left[\begin{matrix}z^n I_{m+1}\\ z^{n-1} I_{m+1}\\ \vdots\\
I_{m+1}\end{matrix}\right] [w^m,w^{m-1},\ldots, 1]^T\label{twovariablechol}.
\end{align}
\end{Lemma}

\begin{proof}
If we substitute the equation
$$
\Phi_{n,m}= \hat L_n(z)[w^m\ \cdots\ 1]^T =\sum_i \hat
L_{n,i}z^i[w^m\ \cdots\ 1]^T
$$
into \eqref{vectorth}, where $\hat L_n(z)$ is an $(m+1)\times
(m+1)$ matrix polynomial of degree $n$, we find, for $j=0,\dots,
n-1$,
\begin{align*}
0&= \left\langle \Phi_{n,m}, \, z^j \left[
\begin{matrix}
w^m\\\vdots \\ 1
\end{matrix}
\right] \right\rangle = \sum^n_{i=0} \hat
L_{n,i} \left\langle z^i \left[
\begin{matrix}
w^m\\\vdots \\ 1
\end{matrix}
\right], \, z^j \left[
\begin{matrix}
w^m\\\vdots \\ 1
\end{matrix}
\right] \right\rangle
\\
& = \sum^n_{i=1} \hat L_{n,i} \left[
\begin{matrix}
\cL_{NM}(z^{i-j}) & \cdots & \cL_{NM}(z^{i-j} w^{-m})
\\
\vdots & & \vdots
\\
\cL_{NM}(z^{i-j} w^m)& \cdots & \cL_{NM}(z^{i-j})
\end{matrix}
\right]
\\
& = \sum^n_{i=1} \hat L_{n,i} \cL_m(z^i,\ z^j) = \cL_m(\hat
L_n(z),\ z^j).
\end{align*}
Similarly,
$$
\langle \Phi_{n,m}, \, \Phi_{n,m}\rangle = I_{m+1} = \cL_m\langle
\hat L_n(z),\, \hat L_n(z)\rangle.
$$
This, coupled with \eqref{leftorth} and the fact that
\eqref{vectpoly} implies that $\hat L_{n,m}$ is upper triangular
with positive diagonal entries, gives \eqref{relmatrix} . Equation \eqref{revmatrix} follows
from \eqref{relmatrix} and \eqref{jsym} while equation \eqref{twovariablechol}  follows
from \eqref{relmatrix} and the definition of $L$.
\end{proof}

Analogous formulas for bivariate orthogonal polynomials in the
reverse lexicographical ordering are obtained by interchanging the
roles on $n$ and $m$.

The function $\cM$ given by \eqref{lfunt} can be used to show that
$\overleftarrow\Phi_{n,m}$ satisfies a minimization condition.
Define $\bar \cM: \prod^{m}_{m+1}\to {\rm Herm}(m+1)$ by
$$
\bar\cM(\Phi)=\langle\Phi^{\dagger}, \, \Phi^{\dagger}\rangle - (\Phi_0+\Phi_0^{\dagger}),
$$
We find

\begin{Lemma}\label{lemmin} The polynomial $\overleftarrow\Phi_{n,m}$ is the unique
minimizer on $\hat\prod^{n,m}_{m+1}$.
\end{Lemma}

\begin{proof} Since $\Phi\in\prod^{n,m}_{m+1}$ can be represented as
$$
\Phi(z,w)=[1,w,...,w^m]\hat\Phi(z)=[1,w,...,w^m]\sum_{i=0}^n\Phi_i z^i,
$$
and from \eqref{relfunct}
\begin{equation*}
\left\langle z^i\left[\begin{matrix} w^m\\\vdots \\ 1
\end{matrix}\right],z^j\left[\begin{matrix} w^m\\\vdots \\ 1
\end{matrix}\right]\right\rangle=\left[\begin{matrix}
\cL_{NM}(z^{i-j}) & \cdots & \cL_{NM}(z^{i-j} w^{-m})
\\
\vdots & & \vdots
\\
\cL_{NM}(z^{i-j} w^m)& \cdots & \cL_{NM}(z^{i-j})
\end{matrix}\right]=\cL_m(z^i,z^j),
\end{equation*}
we find $\bar\cM(\Phi)=\hat\cM(\hat\Phi)$. The  result now follows from
\eqref{revmatrix} and the fact that $R^m_{n,m}\overleftarrow R_n^m(z)$ minimizes
$\hat\cM$ on $\prod^{n}_{m+1}$.

\end{proof}

We can now derive recurrence relations between the various polynomials.

\begin{Theorem}\label{recurrencefor}
Given $\{\Phi_{n,m}\}$ and $\{\Phit_{n,m}\}$, $0\le n\le N$,
$0\le m\le M$, the following recurrence formulas hold
\begin{align}
& A_{n,m}\Phi_{n,m} = z\Phi_{n-1,m} - \hat E_{n,m}\overleftarrow{\Phi}_{n-1,m}^T \label{3.12a}
\\
&\Phi_{n,m}+ A^{\dagger}_{n,m}\hat E_{n,m}(A^T_{n,m})^{-1}\overleftarrow{\Phi}_{n,m}^T=
A^{\dagger}_{n,m}z\Phi_{n-1,m}\label{3.12b}
\\
&\G_{n,m} \Phi_{n,m} = \Phi_{n,m-1} - {\cK}_{n,m} \Phit_{n-1,m},
\label{3.13}
\\
&\G_{n,m}^1 \Phi_{n,m} = w \Phi_{n,m-1} - {\cK}^1_{n,m}\overleftarrow{\Phit}_{n-1,m}^T,
\label{3.14}
\\
&\Phi_{n,m}=I_{n,m} \Phit_{n,m} + \G^{\dagger}_{n,m} \Phi_{n,m-1}
\label{3.15},
\\
&\overleftarrow{\Phi}_{n,m}^T=I_{n,m}^1 \Phit_{n,m} + (\G^1_{n,m})^T\overleftarrow{\Phi}_{n,m-1}^T
\label{3.16},
\end{align}
where
\begin{align}
\hat E_{n,m} & = \langle z\Phi_{n-1,m}, \overleftarrow{\Phi}_{n-1,m}^T\rangle=E_{n,m}J_m=\hat E_{n,m}^T \in
M^{m+1,m+1}, \label{enm}
\\
A_{n,m} & = \langle z\Phi_{n-1,m},\Phi_{n,m}^T\rangle \in
M^{m+1,m+1}, \label{anm}
\\
{\cK}_{n,m} & = \langle \Phi_{n,m-1}, \Phit_{n-1,m}\rangle \in
M^{m,n}, \label{knm}
\\
\G_{n,m} & = \langle \Phi_{n,m-1}, \Phi_{n,m} \rangle \in M^{m,m+1},
\label{gnm}
\\
{\cK}^1_{n,m} & = \langle w \Phi_{n,m-1}, \overleftarrow{\Phit}_{n-1,m}^T\rangle \in M^{m,n},
\label{k1nm}
\\
\G^1_{n,m} & = \langle w \Phi_{n,m-1}, \Phi_{n,m} \rangle \in
M^{m,m+1}, \label{g1nm}
\\
I_{n,m} & = \langle \Phi_{n,m}, \Phit_{n,m}\rangle \in
M^{m+1,n+1}, \label{inm}
\\
I^1_{n,m} & = \langle \overleftarrow{\Phi}_{n,m}^T, \Phit_{n,m}\rangle \in
M^{m+1,n+1}, \label{i1nm}
\end{align}
Similar formulas hold for $\Phit_{n,m}$ and will be denoted by
{\rm(\~{\ref{3.12a}})--(\~{\ref{3.16}})}, etc.
\end{Theorem}

\begin{proof}
Equation \eqref{3.12a} follows from \leref{lemrelmatrix},
equations \eqref{mrecurrence}, \eqref{jsym}, and \eqref{centroahat}\, .Likewise \eqref{3.12b} follows in
an analogous manner from \eqref{frinv}. To prove \eqref{3.13} note
that, because of the linear independence of the entries of
$\Phi_{n,m}$, there is an $m\times(m+1)$ matrix $\Gamma_{n,m}$
such that $\Gamma_{n,m}\Phi_{n,m}-\Phi_{n,m-1}\in
\hat\prod^{n-1,m}_{m}$. Furthermore
$$
\langle \G_{n,m}\Phi_{n,m}-\Phi_{n,m-1},z^iw^j\rangle = 0,
\qquad 0\le i\le n-1\quad 0\le j\le m-1.
$$
Thus \leref{lemvectorthtilde} implies that
$$
\G_{n,m}\Phi_{n,m}-\Phi_{n,m-1}=H_{n,m}\Phit_{n-1,m}.
$$
The remaining recurrence formulas follow in a similar manner.
\end{proof}

\begin{Remark}
As indicated in the proof, formula \eqref{3.12a} follows from the
theory of matrix orthogonal polynomials and so allows us to compute
in the $n$ direction
along a strip of size $m+1$. This formula does not mix the
polynomials in the two orderings. However, to increase $m$ by one
for polynomials constructed in the lexicographical ordering, the
remaining relations show that orthogonal polynomials in the
reverse lexicographical ordering must be used.
\end{Remark}
Using the orthogonality relations from \leref{lemvectorth},
\leref{lemvectorthtilde} and {\eqref{sorthogonal} we find the
following relations.
\begin{Proposition}\label{pr3.6} The following relations hold between the coefficients in
the equations for $\Phit$ and $\Phi$,
\begin{align}
 &{\tilde\cK}_{n,m}={\cK}_{n,m}^{\dagger},\ \tilde I_{n,m}=I_{n,m}^{\dagger},\label{3.24}
\\
& \tilde I^1_{n,m}=(I^1_{n,m})^T, \tilde{\cK}_{n,m}^1=({\cK}^1_{n,m})^T, \label{3.25}
\end{align}
Also
\begin{align}
&A_{n,m}A^{\dagger}_{n,m}=I_m-\hat E_{n,m}\hat E_{n,m}^{\dagger},\label{3.26}
\\
& \G_{n,m}\G_{n,m}^\dagger=I_m-{\cK}_{n,m} {\cK}_{n,m}^\dagger, \label{3.27}
\\
& \G^1_{n,m}(\G^1_{n,m})^\dagger=I_m- {\cK}^1_{n,m} ({\cK}^1_{n,m})^\dagger \label{3.28}
\\
& I_{n,m}I^{\dagger}_{n,m}+\G^{\dagger}_{n,m}\G_{n,m}=I_{m+1} \label{3.29}
\\
& I^1_{n,m}(I^1_{n,m})^\dagger+(\G^1_{n,m})^{\dagger}\G^1_{n,m}=I_{m+1}.
\label{3.30}
\end{align}
\end{Proposition}
\begin{Remark}\label{struct}
The matrix $\G_{n,m}$ has a zero in the entries $(i,j), i\ge j$ and has  positive $(i,i+1)$
entries. Since $\Gamma_{n,m}\Gamma_{n,m}^\dagger=\Gamma_{n,m}U^{\dagger}_m U_m\Gamma_{n,m}^\dagger$
where $U_m$ is the $m\times m+1$ matrix given by\begin{equation}\label{um}
U_m=\left[ \begin{matrix} 0,& I_m \end{matrix}\right],
\end{equation}
we see that $\G_{n,m}U^{\dagger}_m$  is the upper
Cholesky factorization of the right hand side of \eqref{3.27}. From this it is easy to obtain
$\G_{n,m}$. The matrix $\G^1_{n,m}$ has zeroes in the entries $(i,j), i>j$ with positive $(i,i)$
entries. The matrix $I_{n,m}$ has first row and column equal to zero except for a 1 in the $(1,1)$
entry.
\end{Remark}\
The above recurrence formulas also give  pointwise formulas for
the recurrence coefficients. In order to obtain these formulas we
define the $m\times m+1$ matrix $U^1_m$ as
\begin{equation}\label{u1m}
U^1_m=\left[ \begin{matrix}I_m,&0 \end{matrix}\right],
\end{equation}
and the $(n+1)(m+1)\times(n+1)(m+1)$ matrix $P^{n,m}_{rl}$ which
takes monomials in the lexicographical ordering to those in the
reverse lexicographical ordering, i.e.,
\begin{equation}\label{prl}
P^{n,m}_{rl}[z^n w^m,z^n w^{m-1},\ldots ,1]^T=[w^m z^n, w^m z^{n-1},\ldots,1]^T.
\end{equation}
Analogous equations hold for the $n\times (n+1)$ matrices $\tilde
U_n$ and $\tilde U_n^1$.

\begin{Proposition}\label{pr3.8} Let
\begin{equation}\label{eqforphi}
\Phi_{n,m}(z,w)=\Phi_n^m(z)\left[
\begin{matrix} w^m\\\vdots \\ 1
\end{matrix}\right]\ {\rm and}\ \tilde\Phi_{n,m}(z,w)=\tilde\Phi^n_m(w)\left[
\begin{matrix} z^n\\\vdots \\ 1
\end{matrix}\right]
\end{equation}
where
\begin{align}
\nonumber&\Phi_n^m(z) = \Phi^m_{n,n} z^n + \Phi^m_{n,n-1} z^{n-1} +\ldots,\\&\tilde\Phi_m^n(w) =
\tilde\Phi^n_{m,m} w^m + \tilde\Phi^n_{m,m-1} w^{m-1} +\ldots,\label{eqforonephi}
\end{align}
then the following relations hold,
\begin{align}
 &\G_{n,m}=\Phi^{m-1}_{n,n} U_m(\Phi^m_{n,n})^{-1},\label{3.32}
\\
& \G^1_{n,m}=\Phi^{m-1}_{n,n} U_m^1(\Phi_{n,n}^m)^{-1}, \label{3.33}
\\
& {\cK}_{n,m}=-\G_{n,m}I_{n,m}\tilde F_{n,m}, \label{3.34}
\\
& {\cK}^1_{n,m}=-\G^1_{n,m}\bar I^1_{n,m}\bar {\tilde F}^1_{n,m},
\label{3.35}
\\
& I_{n,m}=({\Phi^m_{n,n}}^{\dagger})^{-1}[I_{m+1},0,\ldots,0]C^{-1}_{n,m}{P^{n,m}_{rl}}^T
[I_{n+1},0,\ldots,0]^T (\tilde\Phi^n_{m,m})^{-1},\label{3.36}
\\
& I^1_{n,m}=({\Phi^m_{n,n}}^T)^{-1}[0,\ldots,0,J_{m+1}]C^{-1}_{n,m}{P^{n,m}_{rl}}^T
[I_{n+1},0,\ldots,0]^T (\tilde\Phi^n_{m,m})^{-1},\label{3.37}
\end{align}
where $\tilde F_{n,m}=\tilde\Phi_{m,m}^n U_n^T (\tilde\Phi^{n-1}_{m,m})^{-1}$, and $\tilde
F^1_{n,m}=\tilde\Phi_{m,m}^n (U_n^1)^{T}(\tilde\Phi^{n-1}_{m, m})^{-1}$.
\end{Proposition}
\begin{proof}
Equation \eqref{3.33} follows by equating the coefficients of $z^n$ in \eqref{3.14}
on the left.
The same argument gives \eqref{3.33}.
To show \eqref{3.34} multiply \eqref{3.15} on the left by $\G_{n,m}$
then subtract the resulting equation from \eqref{3.13}.
Now equating the coefficients of $w^m$
gives the result. Equation \eqref{3.35} follows by taking the
transpose of the reverse of equation \eqref{3.14}
then multiplying \eqref{3.16} on the left
by $\bar\G^1_{n,m}$ and subtracting the resulting equations. Equating powers of $w^m$ then
gives the result. Equation \eqref{3.36} follows by equating the highest powers of $w$ in
equation \eqref{3.15} and equation \eqref{3.37} follows in a similar manner from \eqref{3.16}
and the fact that $C_{n,m}$ is a doubly Toeplitz matrix.
\end{proof}

\begin{Remark} From \eqref{relmatrix} and \leref{lemchol} we see that $(\Phi^m_{n,n})^{\dagger}$
is the lower Cholesky factor of
$[I_{m+1},0,\ldots,0]C_{n,m}[I_{m+1},0,\ldots,0]^T $ and a similar
relation holds between $\tilde\Phi^n_{m,m}$and $\tilde C_{n,m}$.
Thus equations \eqref{3.36} and \eqref{3.37} give the relation
between $I_{n,m}$ and $I^1_{n,m}$ and the Fourier coefficients of
$\cL_{N,M}$. These coupled with equations \eqref{3.34} and
\eqref{3.35} relate the Fourier coefficients of $\cL_{N,M}$ to
$\cK_{n,m}$ and $\cK^1_{n,m}$.
\end{Remark}

We now give relations between the coefficients in the recurrence formulas at one level in
terms of those at previous levels.

\begin{Lemma}[Relations for $\cK_{n,m}$]\label{Lemma_K} For $0<n,m$,
\begin{align}
&\G^1_{n,m-1}\cK_{n,m}=\cK_{n,m-1}(\tilde A^{-1}_{n-1,m})^{\dagger}-\cK^1_{n,m-1}\hat{{\tilde E}}_{n-1,m}^{\dagger}(\tilde A^{-1}_{n-1,m})^{\dagger},\label{K1}
\\
&\cK_{n,m}(\tilde\G^1_{n-1,m})^{\dagger}=A^{-1}_{n,m-1}\cK_{n-1,m}-A^{-1}_{n,m-1}\hat E_{n,m-1}\bar\cK^1_{n-1,m}.\label{K2}
\end{align}
\end{Lemma}
\begin{proof}
To show equation \eqref{K1} multiply \eqref{knm} on the left by $\G^1_{n,m-1}$ then use \eqref{3.14} with $m$ reduced by one to obtain,
$$
\G^1_{n,m-1}\cK_{n,m}=\langle w \Phi_{n,m-2}, \Phit_{n-1,m}\rangle.
$$
Eliminating $\Phit_{n-1,m}$ using (\~{\ref{3.12a}}) then applying \eqref{knm} and \eqref{k1nm}
gives \eqref{K1}. Equation \eqref{K2} follows in a analogous manner.
\end{proof}
\begin{Lemma}[Relations for $\cK^1_{n,m}$]\label{Lemma_K1} For $0<n,m$,
\begin{align}
&\G_{n,m-1}\cK^1_{n,m}=\cK^1_{n,m-1}(\tilde A^{-1}_{n-1,m})^T-\cK_{n,m-1}\hat{ \tilde E}_{n-1,m}(\tilde A^{-1}_{n-1,m})^T,\label{K11}
\\
&\cK^1_{n,m}(\tilde\G_{n-1,m})^T=A^{-1}_{n,m-1}\cK^1_{n-1,m}-A^{-1}_{n,m-1}\hat E_{n,m-1}\bar\cK_{n-1,m}.\label{K12}
\end{align}
\end{Lemma}
\begin{proof}
To show \eqref{K11} multiply \eqref{k1nm} on the left by $\G_{n,m-1}$ then use \eqref{3.13} to
obtain
$$
\G_{n,m-1}\cK^1_{n,m}=\langle w \Phi_{n,m-2}, \overleftarrow{\Phit}^T_{n-1,m}\rangle.
$$
Now use (\~{\ref{3.12a}}) with $n$ reduced by one and then equations \eqref{k1nm} and \eqref{knm}
yields \eqref{K11}. Equation \eqref{K12} follows in a similar manner.
\end{proof}
\begin{Lemma}[Relations for $\hat E_{n,m}$]\label{Lemma_E} For $0<n,m$,
\begin{align}
&\G_{n-1,m}\hat E_{n,m}=A_{n,m-1}\cK_{n,m}(I^1_{n-1,m})^{\dagger}+\hat E_{n,m-1}\bar\G^1_{n-1,m},
\label{E1}
\\
&\hat E_{n,m}(\G^1_{n-1,m})^T=I_{n-1,m}(\cK^1_{n,m})^TA^T_{n,m-1}+\G^{\dagger}_{n-1,m}\hat
E_{n,m-1}.\label{E2}
\end{align}
\end{Lemma}
\begin{proof}
To establish \eqref{E1} multiply \eqref{enm} on the left by $\G_{n-1,m}$ then use \eqref{3.13} to
obtain
$$
\G_{n-1,m}\hat E_{n,m}=\langle z\Phi_{n-1,m-1},\overleftarrow\Phi_{n-1,m}^T\rangle.
$$
With the use of \eqref{3.12a} to eliminate $z\Phi_{n-1,m-1}$ we find
$$
\G_{n-1,m}\hat E_{n,m}=A_{n,m-1}\langle \Phi_{n,m-1},\overleftarrow\Phi_{n-1,m}^T\rangle+\hat E_{n,m-1}\langle \overleftarrow\Phi_{n-1,m-1}^T,\overleftarrow\Phi_{n-1,m}^T\rangle.
$$
The second inner product on the right hand side of the above equation evaluates to $\bar \G^1_{n-1,m}$ while the first may be evaluated using \eqref{3.16} followed by \eqref{knm} to give the claimed equation. To obtain \eqref{E2} multiply \eqref{enm} on the right by $(\G^1_{n-1,m})^T$, then use \eqref{3.14} to get
$$
\hat E_{n,m}(\G^1_{n-1,m})^T=\langle z \Phi_{n-1,m},\overleftarrow\Phi_{n-1,m-1}^T\rangle.
$$
Using \eqref{3.12a} to eliminate $\overleftarrow\Phi_{n-1,m-1}^T$ yields
$$
\hat E_{n,m}(\G^1_{n-1,m})^T=\langle z \Phi_{n-1,m},\overleftarrow\Phi_{n,m-1}^T\rangle A^T_{n,m-1}+\langle \Phi_{n-1,m},\Phi_{n-1,m-1}\rangle\hat E^T_{n,m-1}.
$$
Equation \eqref{gnm} can be used to evaluate the second inner product on the right hand side of
the above equation while \eqref{3.16}, (\~{\ref{knm}}) and \eqref{3.24} can be used to obtain
the first inner product.
\end{proof}
\begin{Lemma}[Relation for $\G^1_{n,m}$]\label{Lemma_G} For $0<n,m$,
\begin{align}\label{G11}
\G^1_{n,m}\G_{n,m}^{\dagger}&=I_{n,m-1}\tilde{\hat E}_{n,m}(I^1_{n,m-1})^T+\G_{n,m-1}^{\dagger}\G^1_{n,m-1}\\\nonumber&\quad+\cK^1_{n,m}\bar{\tilde A}_{n-1,m}^{-1}\tilde{\hat E}_{n-1,m}^{\dagger}\tilde A_{n-1,m}\cK_{n,m}^{\dagger}.
\end{align}
\end{Lemma}
\begin{proof}
To show \eqref{G11} multiply \eqref{g1nm} on the left by $\G_{n,m}^{\dagger}$ and use \eqref{3.13} to find,
\begin{equation}\label{g1f}
\G^1_{n,m}\G_{n,m}^{\dagger}=\langle w\Phi_{n,m-1},\Phi_{n,m-1}\rangle -\langle  w \Phi_{n,m-1},\Phit_{n-1,m}\rangle\cK^{\dagger}_{n,m}.
\end{equation}
Eliminating $w\Phi_{n,m-1}$ in the second term on the right hand side of the above equation using \eqref{3.14} then applying (\~{\ref{3.12b}}) gives the third term on the right hand side of \eqref{G11} . In the first term on the right hand side of the above equation substitute the reverse transpose of \eqref{3.16} to find
$$
\langle w\Phi_{n,m-1},\Phi_{n,m-1}\rangle=\langle w\Phi_{n,m-1},\overleftarrow\Phit_{n,m-1}^T\rangle(I^1_{n,m-1})^T+\G^{\dagger}_{n,m-1}\G^1_{n,m-1},
$$
where \eqref{gnm} has been used to obtain the second term on the
right hand side of the above equation. The result may now be
obtained by applying \eqref{3.15} for $\Phi_{n,m-1}$ and then
using (\~{\ref{enm}})
\end{proof}
\begin{Lemma}[Relations for $I_{n,m}$ and $I^1_{n,m}$]\label{Lemma_I}
\begin{align}
&I_{n,m}{\tilde\G}^{\dagger}_{n,m}=-\G^{\dagger}_{n,m}\cK_{n,m},\label{I1}
\\
&I^1_{n,m}=-\bar A^{-1}_{n,m}\hat E^{\dagger}_{n,m}A_{n,m}I_{n,m}+ A^T_{n,m}I^1_{n-1,m}\tilde\G_{n,m},0<n\label{I11}
\end{align}
\end{Lemma}
\begin{proof}
Equation \eqref{I1} follows by multiplying \eqref{inm} on the right by $\tilde\G_{n,m}^{\dagger}$ then using (\~{\ref{3.13}}), \eqref{3.24}, and \eqref{gnm}. In \eqref{i1nm} use (\~{\ref{3.15}}) and \eqref{3.24} to find
$$
I^1_{n,m}=\langle \overleftarrow\Phi_{n,m}^T, \Phi_{n,m}\rangle I_{n,m}+\langle \overleftarrow\Phi_{n,m}^T, \Phit_{n-1,m}\rangle\tilde\G_{n,m}.
$$
The first inner product on the right hand side may be evaluated using \eqref{3.12b}. To evaluate the second inner product eliminate $\overleftarrow\Phi_{n,m}^T$ using the reverse transpose of equation \eqref{3.12b} then use \eqref{i1nm} to obtain the claimed equation.
\end{proof}

\section{Christoffel-Darboux formulas}

The Christoffel-Darboux formula plays an important role in the
theory of one variable scalar and matrix orthogonal polynomials.
Using the connection between two variable orthogonal polynomials
and matrix orthogonal polynomials we derive two variable analogs
of the Christoffel-Darboux formula. These will play an important
role in the theory of two variable stable polynomials discussed
later.

\begin{Lemma} Given $\{\Phi_{n,m}\}$ and $\{\Phit_{n,m}\}$,
\begin{subequations}
\begin{align}
& \overleftarrow\Phi_{n,m}(z,w)\overleftarrow\Phi^{\dagger}_{n,m}(z_1,w_1)-\bar z_1 z\Phi_{n,m}^T(z,w)\Phi^{\dagger}_{n,m}(z_1,w_1)^T
\label{33a}
\\
& \quad =(1-\bar z_1 z)\Phi_{n,m}(z,w)^T\Phi^{\dagger}_{n,m}(z_1,w_1)^T \nonumber
\\
& +\overleftarrow\Phi_{n-1,m}(z,w)\overleftarrow\Phi^{\dagger}_{n-1,m}(z_1,w_1) -\bar z_1 z\Phi_{n-1,m}^T(z,w)\Phi^{\dagger}_{n-1,m}(z_1,w_1)^T
\label{33b}
\\
& \quad =(1-\bar z_1 z)\Phit_{n,m}(z,w)^T\Phit^{\dagger}_{n,m}(z_1,w_1)^T
\nonumber
\\
& +\overleftarrow\Phi_{n,m-1}(z,w)\overleftarrow\Phi_{n,m-1}(z_1,w_1)^T
-\bar z_1 z\Phi_{n,m-1}(z,w)^T\Phi^{\dagger}_{n,m-1}(z_1,w_1)^T.
\label{33c}
\end{align}
\end{subequations}
\end{Lemma}

\begin{proof}
The equality \eqref{33a}$=$\eqref{33b} follows by subtracting
\eqref{mCDeq} with $n$ reduced by one from the original equation then using \leref{lemrelmatrix}. The equality \eqref{33a}$=$\eqref{33c} can be
obtained in the following manner. Let
$$
Z_{n,m}(z,w)=[1,w,\ldots,w^m][I_{m+1},zI_{m+1},\ldots,z^nI_{m+1}],
$$
and $\tilde Z_{n,m}(z,w)$ be given by a similar formula with the
roles of $z$ and $w$, and $n$ and $m$ interchanged. Then from
\leref{lemchol}, \eqref{mCDeq}, and \eqref{relmatrix} we find
\begin{align*}
& \frac{\overleftarrow\Phi_{n,m}(z,w)\overleftarrow\Phi^{\dagger}_{n,m}(z_1,w_1) -
\bar z_1 z\Phi_{n,m}^T(z,w)\Phi^{\dagger}_{n,m}(z_1,w_1)^T}{1-\bar z_1 z}
\\
& \qquad \qquad = Z_{n,m}(z,w) C_{n,m}^{-1}
Z_{n,m}(z_1,w_1)^{\dagger} = \tilde Z_{n,m}(z,w) \tilde C_{n,m}^{-1}
\tilde Z_{n,m}(z_1,w_1)^{\dagger}
\\
& \qquad \qquad = \tilde\Phi_{n,m}^T(z,w)
\tilde\Phi^{\dagger}_{n,m}(z_1,w_1)^T + \tilde Z_{n,m-1}(z,w) \tilde
C_{n,m-1}^{-1} \tilde Z_{n,m-1}(z_1,w_1)^{\dagger}.
\end{align*}
Switching back to the lexicographical ordering in the second term
in the last equation then using Lemma \ref{lemchol} yields the
result.
\end{proof}

As an immediate application of the above lemma we obtain,

\begin{Theorem}[Christoffel-Darboux formula] Given $\{\Phi_{n,m}\}$ and $\{\Phit_{n,m}\}$,
\begin{align*}
& \frac{\overleftarrow\Phi_{n,m}(z,w)\overleftarrow\Phi^{\dagger}_{n,m}(z_1,w_1) -
\bar z_1 z \Phi_{n,m}^T(z,w)\Phi^{\dagger}_{n,m}(z_1,w_1)^T}{1-\bar z_1 z}
\\
& \qquad \qquad = \sum_{k=0}^n \Phi_{k,m}^T(z,w)
\Phi^{\dagger}_{k,m}(z_1,w_1)^T
\\
& \qquad \qquad = \sum_{j=0}^m \Phit_{n,j}^T(z,w)
\Phit^{\dagger}_{n,j}(z_1,w_1)^T.
\end{align*}
\end{Theorem}

In the first line of the above equation the terms $\bar z_1z$
may be replaced by $\bar w_1w$
if we switch to  $\tilde{\Phi}_{n,m}$.

An interesting variant of equation \eqref{33c} is,

\begin{Lemma}\label{le5.1}
\begin{align}
&\Phi_{n,m}(z,w)^T\Phi^{\dagger}_{n,m}(z_1,w_1)^T
-\Phi_{n,m-1}^T(z,w)\Phi^{\dagger}_{n,m-1}(z_1,w_1)^T\nonumber
\\=&\quad\Phit_{n,m}(z,w)^T\Phit^{\dagger}_{n,m}(z_1,w_1)^T
-\Phit_{n-1,m}^T (z,w)\Phit^{\dagger}_{n-1,m}(z_1,w_1)^T.
\label{5.2}
\end{align}
\end{Lemma}

\begin{proof}
Equating the sums in the above Theorem yields
\begin{align}
&\Phi_{n,m}(z,w)^T\Phi^{\dagger}_{n,m}(z_1,w_1)^T-\sum_{j=0}^{m-1} \Phit_{n,j}^T(z,w)\Phit^{\dagger}_{n,j}(z_1,w_1)^T\\& \qquad =\Phit_{n,m}(z,w)^T\Phit^{\dagger}_{n,m}(z_1,w_1)^T-\sum_{j=0}^{n-1} \Phi_{j,m}^T(z,w)\Phi^{\dagger}_{j,m}(z_1,w_1)^T.
\end{align}
Switching to the lexicographical ordering in the sum on the left hand side of the above equation and reverse lexicographical ordering in the sum on the right hand side, extracting the highest terms, then using the Christoffel-Darboux formula to eliminate the remaining sums gives the result.
\end{proof}

\begin{Remark} The above equations can be derived from the
recurrence formulas in the previous sections. However, the
derivation of equation \eqref{33c} is rather tedious.
\end{Remark}

\section{Algorithm}
 In this section we use the relations developed earlier to provide an algorithm that allows us to
 compute the coefficients in the recurrence formula at higher levels in terms of those at
 lower levels plus some indeterminates that are equivalent to the moments. This will allow
 us to construct positive definite doubly Toeplitz matrices. As a byproduct we construct
 the orthogonal polynomials associated with these matrices. More precisely at each level
 we use the new indeterminates and
the coefficients on the levels $(n,m-1)$ and $(n-1,m)$ to
 construct $\cK_{n,m}$ and $\cK^1_{n,m}$.
With this we can construct
 the other coefficients needed to proceed to the next level. The $\hat E_{n,m}$ are closely
 related to the matrix recurrence coefficients needed to compute $C_{n,m}$. Furthermore
 $\Phi_{n,m}$ and $\Phit_{n,m}$ can also be computed. In order to construct the above
 matrices we will have need of the $m\times (m+1)$ matrices $U_m$ and $U^1_m$ given by \eqref{um} and \eqref{u1m} respectively, and the vector
 $e^m_1\in\R^m$ which is the vector with one in the first entry and zeros everywhere else.
From the definition of $\cK^1_{n,m}$ we see that
\begin{align}\label{k1nmr}
{\cK}^1_{n,m} & = \langle w \Phi_{n,m-1},
\overleftarrow{\Phit}_{n-1,m}^T\rangle\nonumber\\& =\langle w\left[\begin{matrix} \phi_{n,m-1}^{m-1}\\ \vdots\\ \phi_{n,m-1}^{0} \end{matrix}\right],\left[\begin{matrix} \overleftarrow{\tilde\phi}_{n-1,m}^{n-1}\\ \vdots\\ z^n\overleftarrow{\tilde\phi}_{n-1,m}^{0} \end{matrix}\right]\rangle\nonumber\\&= c_{-n,-m}d_{n,m}e^m_1(e^n_1)^T+R_{n,m},
\end{align}
where $R_{n,m}$ is an $m\times n$ matrix containing moments $c_{i,j}, \{|i|\le n ,
|j|\le m\}\setminus\{(\pm n,\pm m)\}.$
Likewise with the help of \eqref{eqforphi} and its tilde counterpart we find
\begin{align}\label{knmr}
{\cK}_{n,m}&  = \langle \Phi_{n,m-1}, \Phit_{n-1,m}\rangle\nonumber\\&=c_{-n,m}\Phi^{n,m-1}_{n,m-1}
\left[\begin{matrix} 0\\ \vdots\\1\end{matrix}\right]({\tilde \Phi}^{m,n-1}_{m,n-1}
\left[\begin{matrix} 0\\ \vdots\\1\end{matrix}\right])^{\dagger}+\hat R_{n,m},
\end{align}
where $\hat R_{n,m}$ contains only moments from lower levels.

We proceed as follows, at level $(0,0)$ we have the parameter
$u_{0,0}>0$ which corresponds to $c_{0,0}$. The polynomials
$\Phi_{0,0}$ and $\Phit_{0,0}$ are chosen as $\frac{1}{\sqrt
u_{0,0}}$. From \eqref{inm} and \eqref{i1nm} we see that
$I_{0,0}=1=I^1_{0,0}$. At level $(i,0)$ there is one new parameter
$u_{i,0}$ which can be taken to correspond with the one
dimensional recurrence coefficient i.e. $u_{i,0}=\alpha_i=\hat
E_{i,0}$ corresponding to the $(i,0)$ level and must be less than
one in magnitude. From \eqref{3.26} and the normalization chosen
for the polynomials $A_{i,0}=\sqrt{1-|\hat E_{i,0}|^2}$. This
allows us to compute $\Phi_{i,0}$, and
$\overleftarrow{\Phi}_{i,0}$. The sizes of the matrices given in
\eqref{knm},\eqref{gnm},\eqref{k1nm}, and \eqref{g1nm} show that
$$
\cK_{i,0} =\tilde\cK_{i,0}= \G_{i,0} = \cK^1_{i,0}=\tilde\cK^1_{i,0} = \G^1_{i,0}=0,
$$
where \eqref{3.24}, and \eqref{3.25} have also been used.
Furthermore \eqref{3.29} and \reref{struct} imply that $I_{i,0}=(e^{i+1}_1)^T=\tilde I_{i,0}^{\dagger}$ where \eqref{3.24} has been used. Equation (\~{\ref{3.27}}) implies that $\tilde\G_{i,0}=U_i$ while \eqref{I11} and \eqref{3.25} allow us to compute,
\begin{equation}\label{ii0}
I^1_{i,0}=(\tilde I^1_{i,0})^T=-\left[\begin{matrix}\hat E_{i,0}^{\dagger}& A_{i,0}\hat E_{i-1,0}^{\dagger}&\ldots&-\prod_{j=1}^i A_{j,0}\end{matrix}\right].
\end{equation}
$\Phit_{i,0}$ can now be computed from equation (\~{\ref{3.15}}).


At level $(0,j)$ there is one new parameter $u_{0,j}$ which as above can be taken to correspond with the one dimensional recurrence coefficient i.e. $u_{0,j}=\alpha_j=\tilde{\hat E}_{0,j}$ corresponding to the $(0,j)$ level and must be less than one in magnitude. The analysis for the $(i,0)$ level can be carried over with the roles of the lexicographical and reverse lexicographical orderings interchanged. Thus from (\~{\ref{3.26}}) and the normalization chosen for the polynomials $\tilde A_{0,j}=\sqrt{1-|\tilde{\hat E}_{0,j}|^2}$ which allows us to compute $\Phit_{0,j}$. Again
$$
\tilde\cK_{0,j} =\cK_{0,j}= \tilde\G_{0,j} = \tilde\cK^1_{0,j} =\ \cK^1_{0,j}= \tilde\G^1_{0,j}=0.
$$
Likewise $\tilde I_{0,j}=(e^{j+1}_1)^T=I_{0,j}^{\dagger}$ and $\G_{0,j}=U_j$. Equations (\~{\ref{I11}}) and \eqref{3.25} allow us to compute ${\tilde I}^1_{0,j}$ as above with $i$ and $j$ interchanged as well as the orderings. Equation \eqref{3.15} now allows us to compute $\Phi_{0,j}$.


At level $(n,m)$ with $n,m>0$ there are two new parameters
$u_{n,m}$ and $u_{-n,m}$since $u_{-n,-m}=\overline u_{n,m}$ and
$u_{n,-m}=\overline u_{-n,m}$. These along with the coefficients
on $(n-1,m)$ and $(n,m-1)$ level will be used to compute
$\cK_{n,m}$ and $\cK^1_{n,m}$. This will be sufficient to compute
the remaining coefficients on level $(n,m)$. We begin with,

{\it  Computation of $\cK_{n,m}$.\/} If $n=1, m=1$ then \eqref{knm} shows that
$\cK_{1,1}$ is a scalar which we choose as $\bar u_{1,-1}$. If $m>1$ we see from
\eqref{3.33} and \eqref{u1m} that
$$
\G^1_{n,m-1}\Phi^{m-1}_{n,n}e^m_m=0,
$$
where $e^m_m$ is the m-dimensional vector with zeros in all its entries except the last which is one. Since $\Phi^{m-1}_{n,n}=A_{n,m-1}^{-1}\ldots A_{1,m-1}^{-1}\Phi_{0,0}^{m-1}$ is an upper triangular invertible matrix we find
$$
\G^1_{n,m-1}\Phi^{m-1}_{n,n}((U^1_{m-1})^TU^1_{m-1}+e^m_m(e^m_m)^T)=\G^1_{n,m-1}\Phi^{m-1}_{n,n}(U^1_{m-1})^TU^1_{m-1},
$$
and from \eqref{3.33} $\G^1_{n,m-1}\Phi^{m-1}_{n,n}(U^1_{m-1})^T=\Phi^{m-2}_{n,n}$. Thus equation \eqref{K1} can be written as
\begin{align*}
&U^1_{m-1}(\Phi^{m-1}_{n,n})^{-1}\cK_{n.m}((\tilde{\Phi}^{n-1}_{m,m})^{\dagger})^{-1}\\&=(\Phi^{m-2}_{n,n})^{-1}(\cK_{n,m-1}(\tilde A^{-1}_{n-1,m})^{\dagger}-\cK^1_{n,m-1}\hat{{\tilde E}}_{n-1,m}^{\dagger}(\tilde A^{-1}_{n-1,m})^{\dagger})((\tilde{\Phi}^{n-1}_{m,m})^{\dagger})^{-1}\\&=(\Phi^{m-2}_{n,n})^{-1}(\cK_{n,m-1}-\cK^1_{n,m-1}\hat{{\tilde E}}_{n-1,m}^{\dagger})((\tilde{\Phi}^{n-1}_{m-1,m-1})^{\dagger})^{-1}\\&=H_{n,m-1}.
\end{align*}
In the last equality we have used the fact that $\tilde A_{n-1,m}\tilde{\Phi}^{n-1}_{m,m}=\tilde{\Phi}^{n-1}_{m-1,m-1}$. Likewise,
\begin{align*}
&(\Phi^{m-1}_{n,n})^{-1}\cK_{n,m}(({\tilde\Phi}^{n-1}_{m,m})^{-1})^{\dagger}(U^1_{n-1})^T\\&=(\Phi^{m-1}_{n-1,n-1})^{-1}(\cK_{n-1,m}-\hat E_{n,m-1}\bar\cK^1_{n-1,m})((\tilde{\Phi}^{n-2}_{m,m})^{\dagger})^{-1}\\&={\tilde H}_{n-1,m}.
\end{align*}
If
\begin{equation}\label{defk}
(e^m_m)^T(\Phi^{m-1}_{n,n})^{-1}\cK_{n,m}((\tilde{\Phi}^{n-1}_{m,m})^{\dagger})^{-1}e^n_n=\bar u_{n,-m},
\end{equation}
then $\cK_{n,m}$ can be solved for as,
\begin{align}\label{eqknm}
&\cK_{n,m}\\\nonumber&=\Phi^{m-1}_{n,n}\left(u_{-n,m}e^m_m (e^n_n)^T+(U^1_{m-1})^TH_{n,m-1}+e^m_m (e^m_m)^T\tilde H_{n-1,m}(U^1_{n-1})\right)(\tilde{\Phi}^{n-1}_{m,m})^{\dagger}.
\end{align}
A necessary condition in order to be able to continue is that
$||\cK_{n,m}||<1$.

{\it Computation of $\G_{n,m}$.\/}  Since $\cK_{n,m}$ is presumed to be a
contraction \reref{struct} shows that $\G_{n,m}$ and $\tilde\G_{n,m}$ may be
computed from the upper Cholesky factor of $I-\cK_{n,m}\cK_{n,m}^{\dagger}$ and
$I-\cK_{n,m}^{\dagger}\cK_{n,m}$ respectively.

{\it Computation of $\cK^1_{n,m}$.\/} In $\cK^1_{n,m}$ we see from
\eqref{k1nmr}} that the only new entry is  $(\cK_{n,m})_{1,1}$. If
$n=1,m=1$ set $\cK_{1,1}=\bar u_{1,1}$. If $m>1$ we will show that
all the rows except the first can be obtained from equation
\eqref{K11}. The structure of $\G_{n,m-1}$ implies that
$\G_{n,m-1}e^{m}_1=0$ so that
$$\G_{n,m-1}=\G_{n,m-1}(U_{m-1}^TU_{m-1}+e^m_1(e^m_1)^T)=\G_{n,m-1}U_{m-1}^TU_{m-1} .$$ But $\G_{n,m-1}U_{m-1}^T$ is an invertible matrix, which allows
us to rewrite \eqref{K11} as follows
\begin{equation*}
U_{m-1}\cK^1_{n,m} = (\G_{n,m-1}U_{m-1}^T)^{-1}(\cK^1_{n,m-1}(\tilde A^{-1}_{n-1,m})^T-\cK_{n,m-1}\hat{ \tilde E}_{n-1,m}^T(\tilde A^{-1}_{n-1,m})^T).
\end{equation*}
This gives all the entries in $\cK^1_{n,m}$ except the first row.

Similarly, if $n>1$ we can write
$$
\Gt_{n-1,m} = \Gt_{n-1,m} (U_{n-1}^T U_{n-1} + e^n_1 (e^n_1)^T)= \Gt_{n-1,m} U_{n-1}^T U_{n-1},
$$
i.e., $\Gt_{n-1,m}^T = U_{n-1}^T (U_{n-1}
\Gt^T_{n-1,m})$, and equation \eqref{K12} can be rewritten as
\begin{equation}\label{koneu}
\cK^1_{n,m}U_{n-1}^T =(A^{-1}_{n,m-1}\cK^1_{n-1,m}-A^{-1}_{n,m-1}\hat E_{n,m-1}\bar\cK_{n-1,m})(U_{n-1}\Gt^T_{n-1,m})^{-1}
\end{equation}
Thus the $m\times(n-1)$ matrix $\cK^1_{n,m}U_{n-1}^T$, which is
obtained from $\cK^1_{n,m}$ by deleting the first column, is known
from the previous levels. This allows to compute all entries in
the first row of $\cK^1_{n,m}$ except $(\cK^1_{n,m})_{1,1}$ and we
put,
\begin{equation}\label{konefirstent}
(\cK^1_{n,m})_{1,1}=\bar u_{n,m}
\end{equation}
A necessary condition on the parameters in order to be able to
continue is that $||\cK^1_{n,m}||<1$ which implies that
$|u_{n,m}|<1$

{\it Computation of $\hat E_{n,m}$}

We begin by taking the transpose of \eqref{E2} using the fact that
$\hat E_{n,m}$ is symmetric then multiplying on the left by the
matrix $e^{m+1}_1(e^m_1)^T$. Now multiply \eqref{E1} by $U_m^T$
and add the resulting equations. If $\hat\G_{n-1,m}$ the
$(m+1)\times (m+1)$ matrix obtained by stacking the first row of
$\G^1_{n-1,m}$ on $\G_{n-1,m}$ we find
\begin{align}\label{inve}
\hat\G_{n-1,m}\hat E_{n,m}&=U^T_m(A_{n,m-1}\cK_{n,m}(I^1_{n-1,m})^{\dagger}+\hat E_{n,m-1}\bar\G^1_{n-1,m})\nonumber\\&+e^{m+1}_1(e^m_1)^T(A_{n,m-1}\cK^1_{n,m}I^T_{n-1,m}+\hat E_{n,m-1}\bar{\G}_{n-1,m}).
\end{align}
From the structure of $\G^1$ and $\G$ we see that $\hat\G_{n-1,m}$
is an upper triangular matrix with positive diagonal entries hence
invertible. Thus $\hat E_{n,m}$ can be computed from the above
equation. If $||\hat E_{n,m}||<1$ then $\hat{\tilde E}_{n,m}$ may
be computed from (\~{\ref{E1}}) and (\~{\ref{E2}}). We may also
compute $A_{n,m}$, $\tilde A_{n,m}$ and the polynomials
$\Phi_{n,m}$ and $\tilde{\Phi}_{n,m}$. While the condition that
$\hat E_{n,m}$ be a contraction is necessary and sufficient to be
able to continue it is not optimal in the sense that it does not
take into account the redundancy inherent in the equations giving
$\hat E_{n,m}$. This will be taken into account in the computation
of $\G^1_{n,m}$.

{\it Computation of $\G^1_{n,m}$}. As above we see that \eqref{G11} gives
\begin{align}\label{eqforg1}
&\G^1_{n,m}U^T_m\\\nonumber&=(I_{n,m-1}\hat{\tilde E}_{n,m}(I^1_{n,m-1})^T+\G_{n,m-1}^{\dagger}\G^1_{n,m-1}\\\nonumber&\quad+\cK^1_{n,m}\bar{\tilde A}_{n-1,m}^{-1}\hat{\tilde E}_{n-1,m}^{\dagger}\tilde A_{n-1,m}\cK_{n,m}^{\dagger})(U_m\G^{\dagger}_{n,m})^{-1},
\end{align}
which allows the computation of all the entries of $\G^1_{n,m}$ except the $(1,1)$ entry. Since $(e^m_1)^T I_{n,m-1}= (e^m_1)^T$, $(e^m_1)^T \G^{\dagger}_{n,m-1}= 0$ and likewise with $I_{n,m-1}$ and $\G^{\dagger}_{n,m-1}$ replaced by $\tilde I_{n,m}$ and $\tilde\G_{n,m-1}^T$ respectively we find  with the help of \eqref{inve},
\begin{equation}\label{eg1u}
(e^m_1)^T\G^1_{n,m}U^T_m=(e^m_1)^T H^2_{n,m}+(e^m_1)^T\cK^1_{n,m} H^1_{n,m},
\end{equation}
where
\begin{align}\label{h2nm}
H^2_{n,m}&=I_{n,m-1}((I^1_{n,m-1})^{\dagger}\bar\cK_{n,m}\tilde A^T_{n-1,m}+(\tilde{\G}^1_{n,m-1})^{\dagger} \hat{\tilde E}_{n-1,m})\\\nonumber&\quad\times U_n(\tilde{\hat \G}_{n,m-1}^T)^{-1}(I^1_{n,m-1})^T(U_m\G_{n,m}^{\dagger})^{-1},
\end{align}
and
\begin{align}\label{h1nm}
H^1_{n,m}&=\tilde A^T_{n-1,m}e^n_1(e^{n+1}_1)^T(\tilde{\hat \G}_{n,m-1}^T)^{-1}(I^1_{n,m-1})^T(U_m\G_{n,m}^{\dagger})^{-1}\\\nonumber&\quad+\bar{\tilde A}^{-1}_{n-1,m}\tilde{\hat E}_{n-1,m}^{\dagger}\tilde A_{n-1,m}\cK^{\dagger}_{n,m}(U_m\G_{n,m}^{\dagger})^{-1}.
\end{align}


Thus the first  entry $\G^1$  can be computed using the first row of equation \eqref{eg1u} and \eqref{3.28} which gives
\begin{equation}\label{11g11}
|(\G^1_{n,m})_{(1,1)}|^2=1-(e^m_1)^T H^3_{n,m}(e^m_1),
\end{equation}
where
\begin{equation}\label{h3nm}
H^3_{n,m}=(H^2_{n,m}+\cK^1_{n,m} H^1_{n,m})(H^2_{n,m}+\cK^1_{n,m} H^1_{n,m})^{\dagger}+\cK^1_{n,m}(\cK^1_{n,m})^{\dagger}.
\end{equation}

{\it Computation of the remaining coefficients}. Using the arguments above we see that the relevant part of $I_{n,m}$ may be computed from \eqref{I1} and $I^1_{n,m}$ may be computed from \eqref{I11}. The matrix $\tilde{\G}^1_{n,m} $ can be computed in the same manner as $\G^1_{n,m}$.

\section{Construction of a positive linear functional}

The above algorithm allows us to find a linear functional given
the coefficients in the recurrence formulas. More precisely,

\begin{Theorem}\label{th6.1} Given parameters
$u_{i,j}\in\C$ $0\le i\le n, |j|\le m$,$u_{-i,j}=\bar u_{i,-j}$ we construct
\begin{itemize}
\item  scalars $\hat E_{i,0}$, $i=1,\dots,n$, and
$\tilde{\hat E}_{0,j}$, $j=1,\dots,m$;
\item  $i\times j$ matrices $\cK_{i,j}$
$i=1,\dots,n$, $j=1,\dots,m$;
\item  $i\times j$ numbers $(e^j_1)^T H^3_{i,j}e^j_1$
$i=1,\dots,n$, $j=1,\dots,m$;
\end{itemize}
If
\begin{equation}\label{6.1}
u_{0,0}>0,\ |\hat E_{i,0}|<1, |\tilde{\hat E}_{0,j}|<1,\ ||\cK_{i,j}||<1, \text{ and } {e^j_1}^T H^3_{i,j}e^j_1<1
\end{equation}
then there exists a positive linear functional $\cL$ on $\prod^{n,m}$ such that
\begin{equation}\label{6.2}
\cL(\Phi_{i,m}\Phi_{j,m}^{\dagger})=\delta_{i,j}I_{m+1}
\text{ and }
\cL(\Phit_{n,i}\Phit_{n,j}^{\dagger})=\delta_{i,j}I_{n+1}.
\end{equation}
The conditions \eqref{6.1} are also necessary.
\end{Theorem}

\begin{proof}
We construct the linear functional by induction. First, if $n=m=0$
we set
\begin{equation*}
\cL(1)=u_{0,0}\text{ and }\Phi_{0,0}=\Phit_{0,0}=\frac{1}{\sqrt{u_{0,0}}},
\end{equation*}
and thus $\cL(\Phi_{0,0}\Phi_{0,0}^{\dagger})=\cL(\Phit_{0,0}\Phit_{0,0}^{\dagger})=1$.

If $m=0$, we construct $A_{i,0}=\sqrt{1-|\hat E_{i,0}|^2}$ where
$\hat E_{i,0}=u_{i,0}$. The polynomials $\Phi_{i,0}\ i=0,\ldots n$
are now computed using \eqref{3.12a} and then we define
\begin{equation*}
\cL(\Phi_{i,0}\Phi_{j,0}^{\dagger})=\delta_{i,j}.
\end{equation*}
This gives a well defined positive linear functional on $z^j$ for
$|j|\le n$.

Likewise, if $n=0$, we construct $\Phit_{0,k}$  using
(\~{\ref{3.12a}}) and define
\begin{equation*}
\cL(\Phit_{0,i}\Phit_{0,j}^{\dagger})=\delta_{i,j},
\end{equation*}
which gives the linear functional on $w^j$ for $|j|\le m$.
Thus formula \eqref{6.2} will hold if $m=0$ or $n=0$.

Assume now that the functional $\cL$ is well defined and positive for all levels $0\le i\le n-1$, $0\le j\le m$ and $0\le i\le n$, $0\le j\le m-1$
before $(n,m)$. To ease notation we will use the bracket given in equation~\eqref{innerprod} with $\cL_{N,M}$ replaced by $\cL$. We first extend $\cL$ so that
\begin{equation}\label{6.3}
\langle\Phi_{n,m-1},\Phit_{n-1,m}\rangle=\cK_{n,m}.
\end{equation}
To check that the above equation is consistent with how $\cL$ is
defined on the previous levels, note that from \eqref{K1}
\begin{equation}\label{6.4}
\langle\G^1_{n,m-1}\Phi_{n,m-1},\Phit_{n-1,m}\rangle=\G^1_{n,m-1}\cK_{n,m},
\end{equation}
which follows from the construction of $\cK_{n,m}$ and the
definition of $\cL$ on the previous levels (see \leref{Lemma_K}).
Similarly, using the second defining relation of $\cK_{n,m}$
(i.e., the last row of \eqref{K2}) we see that
\begin{equation}\label{6.5}
\langle\Phi_{n,m-1},\Gt^1_{n-1,m}\Phit_{n-1,m}\rangle=\cK_{n,m}(\Gt^1_{n-1,m})^{\dagger}.\end{equation}
Equations \eqref{6.4} and \eqref{6.5} show that most of \eqref{6.3} is
automatically true. We now define $\cL(z^n w^{-m})$ so that \eqref{defk}
holds which completes \eqref{6.3}.

Using an analogous argument we can use the construction of $\cK^1_{n,m}$ to extend the functional to $z^n w^m$ so that
\begin{equation}\label{6.6}
\cK^1_{n,m}=\langle w\Phi_{n,m-1}, \overleftarrow{\tilde\Phi}^T_{n-1,m}\rangle
\end{equation}
This completes the extension of $\cL$. What remains to show are the equations \eqref{6.2} hold. This is accomplished by first constructing $\tilde{\hat E}_{n,m}$ from (\~{\ref{inve}}). The condition on $(e^m_1)^T H^3_{n,m} e^m_1$ and \eqref{eqforg1} show that the first row of $\G^1_{n,m}$ may be computed and that we may choose
$$
(\G^1_{n,m})_{1,1}>0.
$$
With the first row of $\G^1_{n,m}$ and all of $\G_{n,m}$ (which is calculated from the Cholesky factorization of $\cK_{n,m}\cK_{n,m}^{\dagger}$)\  $\Phi_{n,m}$ may be constructed from \eqref{3.13} and \eqref{3.14}.
Equations \eqref{6.3} and \eqref{6.6} coupled with \eqref{3.13},
\eqref{3.14} and the orthogonality relations on the previous levels show that
\begin{equation*}
\langle\G_{nm}\Phi_{n,m},\ \tilde\Phi_{n-1,k}\rangle=0,
 \quad k=0,1,\dots,m,
\end{equation*}
and
\begin{equation}\label{orthpl}
\langle(e^m_1)^T\Gamma^1_{nm}\Phi_{n,m},\ w^k
\left[\begin{matrix}z^{n-1}\\ \vdots\\1\end{matrix}\right]\rangle
=0,\quad k=1,2,\dots, m .\end{equation} Equations \eqref{orthpl}
and \eqref{3.14} show
$$
0=\langle(e^m_1)^T\Gamma^1_{n,m}\Phi_{n,m},\ \overleftarrow{\tilde \Phi}^T_{n-1,m}\rangle
=\langle(e^m_1)^T\Gamma^1_{nm}\Phi_{n,m},\ \left[
\begin{matrix} z^{n-1}\\ \vdots\\1\end{matrix}\right]\rangle.$$
The fact that $\overleftarrow{\tilde \Phi}^T_{n-1,m}$ has an invertible coefficient
multiplying $\left[\begin{matrix} z^{n-1}\\ \vdots\\1\end{matrix}\right]$
has been used to obtain the second equality in the above equation.  The
above implies
$$\langle\Phi_{n,m},\ \tilde\Phi_{n-1,k}\rangle=0, \qquad k=0,1,\dots, m, $$
which in turn implies that
$$\langle \Phi_{n,m},\ \Phi_{j,m}\rangle=0, \qquad j=0,1,\dots, n-1 .$$
To show that
$$\langle\Phi_{n,m},\ \Phi_{n,m}\rangle = I_{m+1}$$
we note that equations \eqref{3.13}, \eqref{3.14} and \eqref{G11}
imply
$$\langle(e^m_1)^T\Gamma^1_{nm}\Phi_{n,m},\ \Gamma_{nm}\Phi_{n,m}\rangle=
(e^m_1)^T\Gamma^1_{n,m}\Gamma^{\dagger}_{n,m}$$
and \eqref{3.28} implies
$$
\langle(e^m_1)^T\Gamma^1_{nm}\Phi_{n,m},\ (e^m_1)^T\Gamma^1_{nm}\Phi_{nm}\rangle
 = (e^m_1)^T\Gamma^1_{nm}(\Gamma^1_{nm})^{\dagger}(e^m_1).
$$
Thus $\cL$ is a positive linear functional. The  orthogonality relations for the polynomials $\Phit_{i,j}$  now follow.
\end{proof}

Let $C(\TT^2)$ denote the set of continuous functions on the bi-circle, above Theorem now allows,
\begin{Theorem}\label{th6.2} Given parameters
$u_{i,j}\in\C$ with $u_{i,-j}=\bar u_{-i,j}$. If equations \eqref{6.1} hold for
all $0\le i,j$ then there exists a positive measure $\mu$ supported on the bi-circle such that for any $f\in C(\TT^2)$,
$$\cL(f)=(\frac{1}{2\pi})^2\int_{\TT^2} f(\theta,\phi) d\mu(\theta,\phi)$$
\end{Theorem}

\begin{proof} From the  hypotheses imposed above, Theorem \ref{th6.1}   shows that $C_{n,m}$ is positive definite for all $n$ and $m$ so the result  follows from
Bochner's Theorem \cite[Section 1.4.3]{R}.
\end{proof}

\begin{Remark}
The above construction gives a criteria for the existence of
a one step extension of the functional. That is, given moments so
that there exists a positive linear functional on $\prod^{n-1,m}
\cup \prod^{n,m-1}$, any set
$$
\{ u_{n,m}, \, u_{-n,m}\}\ u_{-n,-m}=\bar u_{n,m}\quad u_{n,-m}=\bar u_{-n,m}
$$
that satisfies \eqref{6.1} can be used to extend the functional to
$\prod^{n,m}$. However it is not difficult to construct examples
where no extension exists. See section 8.2
\end{Remark}

\section{Two Variable Stable polynomials and  Fej\'er-Riesz factorization}

In this section we study the consequences of $\cK_{n,m}=0$. This will make a connection with the
results in \cite{GW1} on stable polynomials and the Fej\'er-Riesz factorization theorem.
We say that a polynomial $p(z,w)$ is stable if $p(z,w)\ne 0, \ |z|\le 1,\ |w|\le1$.
A polynomial $p$ is of degree $(n,m)$ if
$$
p(z,w)=\sum_{i=0}^n \sum_{j=0}^m k_{i,j} z^i w^j,
$$
with $k_{n,m}\ne0$. Finally we say that the polynomial $p_{n,m}$
of degree $(n,m)$ has the spectral matching property (up to
$(n,m)$) if
$$
\cL(z^k w^j)=\frac{1}{(2\pi)^2}\int_{\TT^2}\frac{z^k
w^j}{|p_{n,m}(z,w)|^2}d\theta d\phi,\quad  z=e^{i\theta}\
,w=e^{i\phi},
$$
for $|k|\le n$, $|j|\le m$.

\begin{Lemma}\label{christdarblike} Suppose that $\cL$ is a positive definite linear functional on
$\prod^{n,m}$ and $\cK_{n,m}=0$, then,
\begin{align}
& \overleftarrow\phi^m_{n,m}(z,w)\overline{\overleftarrow\phi^m_{n,m}(z_1,w_1)}-\phi_{n,m}^m(z,w)\overline{\phi^m_{n,m}(z_1,w_1)}\nonumber\\ & \quad=(1-w\bar w_1)\overleftarrow\Phi_{n,m-1}(z,w)\overleftarrow\Phi^{\dagger}_{n,m-1}(z_1,w_1)\nonumber\\ &
\quad +(1-z\bar z_1)\Phit_{n-1,m}(z,w)^T\Phit^{\dagger}_{n-1,m}(z_1,w_1). \label{chlike}
\end{align}
\end{Lemma}

\begin{proof} If $\cK_{n,m}=0$ then \eqref{3.13} shows that $\G_{n,m}\Phi_{n,m}(z,w)=\Phi_{n,m-1}(z,w)$. Thus $\overleftarrow\Phi_{n,m}(z,w)\G^{\dagger}_{n,m}=w\overleftarrow\Phi_{n,m-1}(z,w)$. Also \eqref{3.27} implies that $(\G_{n,m})_{(i,i+1)}=1,\ i=1..m$ with all other entries zero. Thus we find
\begin{align}
&\overleftarrow\Phi_{n,m}(z,w)\overleftarrow\Phi_{n,m}(z_1,w_1)^{\dagger}\nonumber\\ &=\overleftarrow\phi^m_{n,m}(z,w)\overline{\overleftarrow\phi^m_{n,m}(z_1,w_1)}+\overleftarrow\Phi_{n,m}(z,w)\G_{n,m}^{\dagger}\G_{n,m}\overleftarrow\Phi_{n,m}(z_1,w_1)^{\dagger}\nonumber\\ &= \overleftarrow\phi^m_{n,m}(z,w)\overline{\overleftarrow\phi^m_{n,m}(z_1,w_1)}+w\bar w_1\overleftarrow\Phi_{n,m-1}(z,w)\overleftarrow\Phi_{n,m-1}(z_1,w_1)^{\dagger}.
\end{align}
From \eqref{33c} we find
\begin{align*}
& \overleftarrow\phi^m_{n,m}(z,w)\overline{\overleftarrow\phi^m_{n,m}(z_1,w_1)}-z\bar z_1\phi_{n,m}^m(z,w)\overline{\phi^m_{n,m}(z_1,w_1)} \\ & \quad=(1-w\bar w_1)\overleftarrow\Phi_{n,m-1}(z,w)\overleftarrow\Phi^{\dagger}_{n,m-1}(z_1,w_1)\\ &
\quad +(1-z\bar z_1)\Phit_{n,m}(z,w)^T\Phit^{\dagger}_{n,m}(z_1,w_1).
\end{align*}
Using (\~{\ref{3.13}}) and the fact that $\tilde\phi^n_{n,m}(z,w)=\phi^m_{n,m}$ gives the result.
\end{proof}
We now have
\begin{Theorem}\label{stablethm} Suppose that $\cL$ is a positive definite linear functional on
$\prod^{n,m}$ and $\cK_{n,m}=0$ then $\overleftarrow\phi^m_{n,m}(z,w)$ is stable and,
$$\cL(e^{-ik\theta}e^{-il\phi})=(\frac{1}{2\pi})^2\int_{\TT^2}\frac{e^{-ik\theta}e^{-il\phi}}{|\phi^m_{n,m}(e^{i\theta},e^{i\phi})|^2}d\theta d\phi,\quad |k|\le n,\ |l|\le m.
$$
Conversely if $\pi_{n,m}(z,w)$ is a polynomial of degree $(n,m)$ such that $\overleftarrow\pi_{n,m}$ is stable and
$$\cL(e^{-ik\theta}e^{-il\phi})=(\frac{1}{2\pi})^2\int_{\TT^2}\frac{e^{-ik\theta}e^{-il\phi}}{|\pi_{n,m}(e^{i\theta},e^{i\phi})|^2}d\theta d\phi,\quad |k|\le n,\ |l|\le m,
$$
then $\cK_{n,m}=0$.
\end{Theorem}

\begin{proof} If $\cL$ is positive definite and $\cK_{n,m}$ is equal to zero
then \leref{christdarblike} shows that $\overleftarrow
\phi^n_{n,m}(z,w)$ satisfies \eqref{chlike}. The first part of the
result now follows from the proof of Theorem 2.3.1 in \cite{GW1}.
To show the second part let
$f(z,w)=1/|\overleftarrow\pi_{n,m}(z,w)|^2, |z|=1=|w|$ be the
spectral density function associated with $\pi_{n,m}$. Then from
equation~2.1.5 in \cite{GW1} and \leref{lemrelmatrix} we find that
$$\overleftarrow\Phi_{n,m}(z,w)=[\overleftarrow\pi_{n,m}(z,w),
w\overleftarrow\Phi_{n,m-1}(z,w)].$$
But this implies that $\Phi_{n,m}(z,w)=[\pi_{n,m}(z,w),\Phi_{n,m-1}(z,w)^T]^T$. Hence from \eqref{3.13} $\cK_{n,m}=0$.
\end{proof}

This leads to the following alternative proof of the two-variable
Fej\'er-Riesz Theorem in \cite{GW1}.

\begin{Theorem}\label{twovarfr}
 Suppose that $f(z,w)=\sum^n_{k=-n}\sum^m_{l=-m}f_{kl}z^kw^l$ is positive for $|z|=|w|=1$.  Then there exists a polynomial
 $$p(z,w)=\sum^n_{k=0}\sum^m_{l=0}p_{kl}z^kw^l$$ with $p(z,w)\ne 0$ for
 $|z|,|w| \le 1$, and $f(z,w)=|p(z,w)|^2$ if and only if $\cK_{n,m}=0$.
\end{Theorem}

\begin{proof} For $g\in C(\TT^2)$ let $\cL(g)=\frac{1}{(2\pi)^2}\int_{\TT^2}\frac{g(\theta,\phi)}{|p(e^{i\theta},e^{i\phi})|^2}d\theta d\phi$. Then $\cL$ is a positive definite linear functional on $\TT^2$.  The necessary part of the above result now follows from \thref{stablethm}. The sufficiency also follows from the above Theorem and the maximal entropy condition \cite{BN}.
\end{proof}

An alternative approach for finding a factorization as above may
be done using the notion of intersecting zeros (see \cite{GW2}).
Also, the question of factorizing a nonnegative trigonometric
polynomial as a modulus square of an outer polynomial was
addressed in \cite{DW}, allowing for generalizations in the
operator valued case. When such a factorization of the desired
degree does not exist one can approximate the trigonometric
polynomial with one that does have the desired factorization. This
question was pursued in \cite{HW}.

The vanishing of $\cK_{n,m}$ has the following geometric
interpretation.

\begin{Lemma}\label{geom} Suppose $\cL$ is positive definite on $\prod^{n,m}$ then
$\cK_{n,m}=0$ if and only if for $\Phi_{n,m-1}$ constructed as in \eqref{vectpoly} ,
\begin{equation}\label{geomeq}
\langle\Phi_{n,m-1}, z^i w^j\rangle=0\quad j=m\ ,0\le i\le n-1 .
\end{equation}
\end{Lemma}

\begin{proof} The definition of  $\Phi_{n,m-1}$ shows that it is already orthogonal
to $z^i w^j$, $0\le i <n,\ 0\le j \le m-1$.
The remaining orthogonality conditions show that $\Phi_{n,m-1}$ is orthogonal to all the
monomials in $\tilde\Phi_{n-1,m}$. Thus the sufficiency part of the Theorem follows from
\eqref{knm}. To see the necessary part note that from the definition of $\Phi_{n,m-1}$
\begin{equation}
\cK_{n,m}=\langle \Phi_{n,m-1}, \tilde\Phi^{n-1}_{n-1,m}\left[
\begin{matrix} z^{n-1}\\\vdots \\ 1\end{matrix}\right] w^m\rangle,
\end{equation}
with $\tilde\Phi^{n-1}_{n-1,m}$ an invertible matrix. Thus \eqref{geomeq} follows.
\end{proof}

Unfortunately at this point we are unable to see what the
condition $\cK_{n,m}=0$ implies for $u_{i,j},\ |i|\le n,\ |j|\le
m$ except for $u_{n,-m}=0$ which follows from
equation~\eqref{defk}. We can however get a partial
characterization for when a positive measure on the bi-circle can
be written as the reciprocal of the magnitude square of a stable
polynomial. We begin with the following auxiliary result.

\begin{Lemma}\label{eij} If $\hat E_{i,j}=0$ then the first column of $K^1_{i,j}$ is equal to zero, in particular $u_{i,j}=0$. If $\hat E_{i,j},\ \cK_{i,j}$ and $\cK_{i-1,j}$ are zero then so is $\hat E_{i,j-1}$. Conversely if $\cK_{i,j},\ K_{i-1,j}$, $\hat E_{i,j-1}$ and $u_{i,j}$ are zero then $\hat E_{i,j}=0$. In both cases $\cK^1_{i,j}=[0,\cK^1_{i-1,j}]$. Likewise if  $\hat {\tilde E}_{i,j}=0$ then the first row of $K^1_{i,j}$ is equal to zero. If $\hat{\tilde E}_{i,j},\ \cK_{i,j}$ and $\cK_{i,j-1}$ are zero then so is $\hat{\tilde E}_{i,j-1}$. Conversely if $\cK_{i,j},\ K_{i,j-1}$, $\hat{\tilde E}_{i-1,j}$ and $u_{i,j}$ are zero then $\hat{\tilde E}_{i,j}=0$. In both cases $\cK^1_{i,j}=[0,(\cK^1_{i,j-1})^T]^T$.
\end{Lemma}

\begin{proof} If $\hat E_{i,j}=0$ then equation~\eqref{E2}, and \reref{struct} show that the first column of $\cK^1_{i,j}$ is zero. If $\cK_{i-1,j}$ is equal to zero then \eqref{I1} shows that all the entries of $I_{i-1,j}$ are all zero except for a one in the first entry. Thus \eqref{E1} and \eqref{E2} imply that if $\hat E_{i,j}=0$, $\cK_{i,j}=0$, and $\cK_{i-1,j}=0$ then $\hat E_{i,j-1}\bar\G^1_{i-1,j}=0$ and $\hat E_{i,j-1}\bar \G_{i-1,j}=0$. Following the argument in the construction of $\hat E_{n,m}$ we see that $\hat E_{i,j-1}=0$. The above hypothesis on $\cK_{i-1,j}$ shows that $\tilde\G_{i-1,j}=U_{i-1}$ thus \eqref{koneu} and the fact that the first column of $\cK^1_{i,j}$ is zero gives $\cK^1_{i,j}=[0,\cK^1_{i-1,j}]$. The converse statement follows from equation~\eqref{inve} . The remaining statements follow in an analogous fashion using \prref{pr3.6} .
\end{proof}

\begin{Lemma}\label{abscont} Let $\mu$ be a positive measure on the bi-circle. Then $\mu$ is purely absolutely continuous with respect to Lebesgue measure and
$$d\mu(\theta,\phi)=\frac{1}{|p_{n,m}|^2}d\theta d\phi,
$$
where $p_{n,m}$ is a polynomial of degree $(n,m)$ with $\overleftarrow p_{n,m}(z,w)$ stable if and only if
$$
\cK_{i,j}=0,\hat E_{i+1,j}=0, {\rm and}\ \hat{\tilde E}_{n, j+1}=0, i\ge n, j\ge m.
$$
\end{Lemma}

\begin{proof} Suppose that $d\mu = \frac{1}{|p_{n,m}(z,w)|^2} d\theta d\phi$ with $\overleftarrow p_{n,m}$ stable, then the sequence
$$\{\psi_{i,j}(z,w)\}\quad \psi_{i,j}(z,w)=z^{i-n} w^{j-m} p_{n,m}(z,w),\ i\ge n,\ j\ge m,
$$
is a set of polynomials with degrees $(i,j)$ respectively such
that $\overleftarrow \psi_{i,j}=\overleftarrow p_{n,m}$ are stable
and have the spectral matching property. Thus \thref{stablethm}
implies that $\cK_{i,j}=0$ for $i\ge n,\ j\ge m$. Since
$\overleftarrow\psi_{i+1,j}=\overleftarrow \psi_{i,j}$, $i\ge n,\
j\ge m$ we see from  equation~2.1.5 in \cite{GW1} that
$\overleftarrow\Phi_{i+1,j}=\overleftarrow\Phi_{i,j}$ for $i\ge n\
, j\ge m$. This implies that $A_{i+1,j}= I_{j+1}$ so that
equation~\eqref{3.12a} shows that $\hat E_{i+1,j}=0,\ i\ge n, j\ge
m$. Since
$\overleftarrow{\tilde\phi}^i_{i,j+1}=\overleftarrow{\tilde\phi}^i_{i,j}$,
$i\ge n$, $j\ge m$ the preceding argument shows that $\tilde
E_{i,j+1}=0,\ i\ge n, j\ge m$. This proves the necessary part.

To prove sufficiency note that if $\cK_{i,j}=0, i\ge n,\ j\ge m$
there exist polynomials $\psi_{i,j}$ of degree $(i,j)$ where
$\overleftarrow\psi_{i,j}$ is a stable polynomial which has the
spectral matching property. In order to show that
$\overleftarrow\psi_{i,j}=\overleftarrow\psi_{n,m}$ we note that
since $\hat E_{i+1,j}=0$ equation~\eqref{3.12a} implies that
$\Phi_{i+1,j}=\Phi_{i,j}$, $i\ge n,\ j\ge m$. Furthermore
$\hat{\tilde E}_{n, j+1}=0,\ j\ge m$ implies that
$\tilde\Phi_{n,j+1}=\tilde\Phi_{n,m}$. Since
$\psi_{i,j}=\phi^j_{i,j}=\tilde\phi^i_{i,j}$ for $i\ge n,\ j\ge m$
the result follows. \end{proof}

The conditions on $\cK_{i,j}$, $E_{i,j}$ and $\tilde E_{n,j}$
given in Lemma \ref{abscont} are not optimal since they are
redundant. Some of this redundancy is removed in the next theorem.

\begin{Theorem}\label{verbtwod} Let $\mu$ be a positive measure on the bi-circle. Then $\mu$ is purely absolutely continuous with respect to Lebesgue measure and $d\mu=\frac{d\theta d\phi}{|p_{n,m}|^2}$ where $p_{n,m}$ is a polynomial of degree $(n,m)$ with $\overleftarrow p_{n,m}$ stable if and only if

\item{a.} $\cK_{n,j}=0$, $\tilde{\hat E}_{n-1,j+1}=0$ and $u_{n,j+1}=0$, $j\ge m$,
\item{b.} $\cK_{i,m}=0$, $\hat E_{i,m-1}=0$ and $u_{i,m}=0$, $i>n$
\item{c.} $u_{|i|,j}=0,i>n,\ j>m$.
\end{Theorem}

\begin{Remark} Equation~\eqref{defk} and \leref{eij} show that $u_{-n,j}$, $u_{-i,m}$, $u_{n-1,j+1}$ and $u_{i,m-1}$ are also equal to zero for $j\ge m, i>n$.
\end{Remark}

\begin{proof} If $\mu$ has the form indicated in the hypotheses then \leref{abscont} says that $\cK_{i,j}=0, \ i\ge n,\ j\ge m$ which coupled with equation~\eqref{defk} implies that $u_{-i,j}=0,\ i\ge n ,\ j\ge m$ the remaining conditions on the coefficients follow from \leref{eij}. If the coefficients obey a.-c. then
\leref{eij} shows that $\hat E_{i+1,m}$ and $\hat{\tilde E}_{n,j+1}$, are equal to zero for $i\ge n$ and $j\ge m$. Since $\hat E_{n,m+1}=0$, $\hat{\tilde E}_{n+1,m}=0$ and by hypothesis $u_{-n-1,m+1}=0$, equation~\eqref{eqknm} shows that $\cK_{n+1,m+1}=0$. With this \leref{eij} shows that $\hat E_{n+1,m+1}=0$ and $\tilde{\hat E}_{n+1,m+1}=0$. The result now follows by induction.
\end{proof}

It is possible to modify slightly the hypotheses of Theorem
\ref{verbtwod} to obtain a statement just on the coefficients in
the recurrence formulas.

\begin{Theorem}\label{param} Suppose $u_{i,j}$ are given so that equations \eqref{6.1}
are satisfied for $0\le i\le n\ , |j|\le m$ and so that a-c in
Theorem \ref{verbtwod} hold.  Then for $f\in C(\TT^2)$,
$$
\cL(f)=(\frac{1}{2\pi})^2\int_{\TT^2} f(\theta,\phi) d\mu(\theta,\phi),
$$
where $\mu$ is absolutely continuous with respect to Lebesgue
measure with $d\mu=\frac{d\theta d\phi}{|p_{n,m}|^2}$. Here
$p_{n,m}$ is a polynomial of degree $(n,m)$ with $\overleftarrow
p_{n,m}$ stable.
\end{Theorem}

\begin{proof} From \thref{th6.1} that there exists a positive definite linear functional on $\prod^{n,m}$ with the above parameters and from \thref{stablethm} the functional has the representation
$$\cL(e^{-ik\theta}e^{-il\phi})=(\frac{1}{2\pi})^2\int_{\TT^2}\frac{e^{-ik\theta}e^{-il\phi}}{|p_{n,m}(e^{i\theta}e^{i\phi})|^2}d\theta d\phi,\quad |k|\le n,\ |l|\le m.
$$
with $p_{n,m}$ a polynomial of degree $(n,m)$ with $\overleftarrow p_{n,m}$ stable. The result now follows from \thref{verbtwod}.
\end{proof}

\section{Examples}

We now give some examples that illustrate various aspects of the results presented earlier. We begin with the case $n=1, m=1$ with $u_{0,0}=1$, $\cK_{1,1}=u_{-1,1}=0$ and $\cK^1_{1,1}=\bar u_{1,1}$.
From \thref{th6.1} we see that we must chose $|u_{0,1}|<1$ and $|u_{1,0}|<1$. Since $\cK_{1,1}=0$ the only remaining condition for $\cL$ to be a positive linear functional on $\prod^{1,1}$ is for ${e^1_1}^T \tilde H^3_{1,1}e^1_1<1$. From equation~\eqref{eqforg1} we see that $\G^1_{1,1}U^T_1=I_{1,0}\hat{\tilde E}_{1,1}(I^1_{1,0})^T$. The construction of $I_{1,0}$, $I^1_{1,0}$ and equation \eqref{inve} shows that ${e^1_1}^T \tilde H^3_{1,1}e^1_1<1$ is given by
$$
a|u_{1,1}|^2+ b(\bar u_{1,1}\bar u_{0,1} u_{1,0}+ u_{1,1}u_{0,1}\bar u_{1,0})+c <1,
$$
with $a=\frac{1-|u_{0,1}u_{1,0}|^2}{1-|u_{1,0}|^2}$, $b=\frac{\sqrt{1-|u_{0,1}|^2}}{\sqrt{1-|u_{1,0}|^2}}$ and $c=|u_{0,1}|^2$.
This simplifies to
$$
|\hat u_{1,1}|<1,
$$
where
$$
\hat u_{1,1}= \frac{(1-|u_{0,1} u_{1,0}|^2)u_{1,1}}{\sqrt{1-|u_{0,1}|^2}\sqrt{1-|u_{1,0}|^2}}+u_{0,1}\bar u_{1,0}.
$$
Thus from \thref{th6.1} and  \thref{stablethm} we see that with $u_{0,0}=1$,
$$\cL( e^{-ik\theta} e^{-ij\phi})=(\frac{1}{2\pi})^2\int_{\TT^2}\frac{e^{-ik\theta} e^{-ij\phi}}{|\phi_{1,1}(e^{i\theta},e^{i\phi})|^2}d\theta d\phi\ |k|\le1,\ |j|\le1,
$$ where $\phi_{1,1}$ constructed using equation \eqref{3.13} and the top row of \eqref{3.14} is a polynomial of degree (1,1) with $\overleftarrow\phi_{1,1}$ stable
if and only if $|u_{0,1}|<1$, $|u_{1,0}|<1$, $u_{-1,1}=0$ and $|\hat u_{1,1}|<1$.
Furthermore if we set $u_{j,0}$, $u_{0,j}$, $u_{i,j}$ equal to zero for $i>1,|j|>1$ then \thref{param} shows that the above representation for $\cL$ extends to all continuous functions on $\TT^2$.

We can also use the previous results to investigate contractive Toeplitz matrices. In this case we find
\begin{equation}
C_{1,1}= \left[\begin{matrix} I & C_{-1}\\
C_1&I
\end{matrix}
\right] ,
\end{equation}
where $C_{-1}=C_1^{\dagger}$ is a $2\times 2$ Toeplitz matrix. In this case $u_{0,0}=1$ and $u_{0,1}=0$ so that $\tilde
E_{0,1}=0$, $\tilde A_{0,1}=1$. Since $\cK_{1,1}=u_{-1,1}$, we find $\G_{1,1}=[0,\sqrt{1-|u_{-1,1}|^2}]$. This plus the computation of $\tilde\G^1_{1,0}$ described in the construction of $\cL$
yields
\begin{align}
&I=(e^1_1)^T\tilde H^3_{1,1}(\tilde H^3_{1,1})^{\dagger}(e^1_1)\\\nonumber&=(1+d)|u_{1,1}|^2 + d(u_{1,1}u_{-1,1}+\bar u_{1,1}\bar u_{-1,1})+d|u_{-1,1}|^2<1,\end{align}
where
$$d=\frac{|u_{1,0}|^2}{(1-|u_{1,0}|^2)(1-|u_{-1,1}|^2)} .
$$
By completing the square this can be simplified to
$$|\hat u_{1,1}|<1$$
where
$$\hat u_{1,1}=(1+d)\sqrt{1-|u_{1,0}|^2}u_{1,1}+d_1\bar u_{-1,1},$$
and
$$d_1=\frac{|u_{1,0}|^2}{(1-|u_{-1,1}|^2)\sqrt{1-|u_{1,0}|^2}}.$$
which puts constraints on $u_{-1,1}$.
Thus we find the conditions for $\cL$ to be a positive linear functional and hence $C_1$ to be a contractive Toeplitz matrix are $|u_{1,0}|<1$, $|u_{-1,1}|<1$ and $|\hat u_{1,1}|<1$.
These constraints may not be strong enough to allow $\cL$ to be extended. To see this suppose $n=1, m=2$, $u_{0,2}=0$ and $u_{1,0}=0$.
It is not difficult to see then that $\hat E_{1,1}={\rm diag}(\bar u_{1,1}, u_{-1,1})$. With $u_{1,0}=0$ the constraint on $\hat u_{1,1}$ above reduces to $|u_{1,1}|<1$. However
$$
K_{1,2}=\left(\begin{matrix}\frac{u_{-1,1}}{(1-|u_{1,1}|^2)^{1/2}}\\ \frac{u_{-1,2}}{(1-|u_{-1,1}|^2)^{1/2}}\end{matrix}\right),
$$
so we see that in order for $K_{1,2}$ to be a contraction $\frac{|u_{-1,1}|}{\sqrt{1-|u_{1,1}|^2}}<1$, which may not be satisfied.

\end{document}